# Ancient Computers

## Part I - Rediscovery

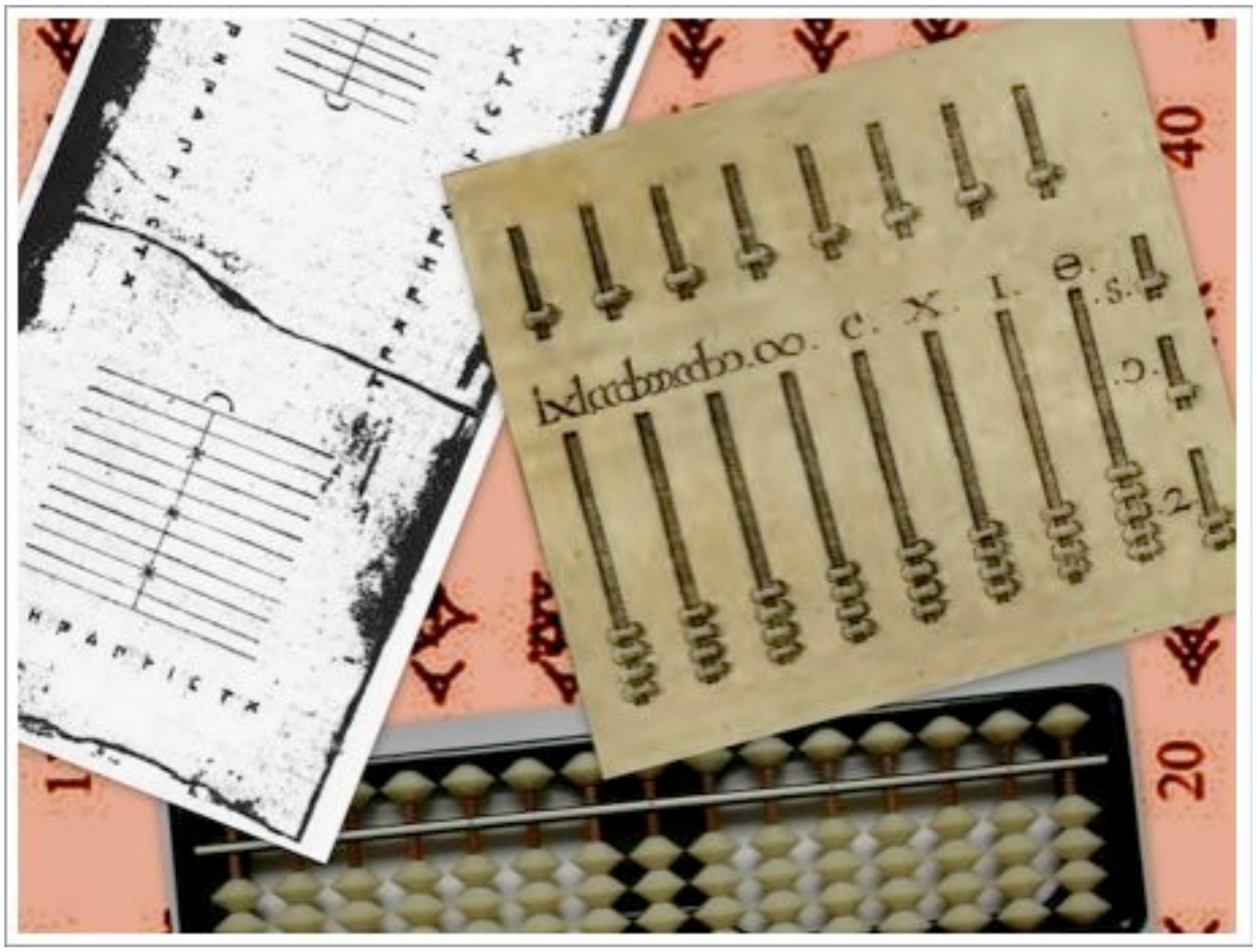

Stephen Kent Stephenson

## Table of Contents

## Author's Forward

In Tokyo in 1964 I bought a Soroban with Kojima's book "The Japanese Abacus: Its Use and Theory"; an event that sparked my interest in abaci ... and in computers.

After getting an M.Eng.(Elect.) degree at Cornell University[1], my 30 year career included working on the design and construction of nuclear power plants, missile systems software engineering, and industrial and engineering computer systems sales and systems engineering.

Deciding to become a high school math teacher at the end of 2000, I took a History of Math course as part of an M.Ed. Program at UMassLowell.[2] I was struck by how easy it would be to use ancient Roman, Greek, Egyptian, and Babylonian numerals to record abaci calculation results. Prof. Gonzalez asked, "Yes, but how would you do multiplication and division?"

So as a hobby, I've worked the last 10 years to (re)discover the schematics and programming rules of the computers the Ancients used to do their accounting and engineering to support and empower the greatest empires in human history.

I hope you find Ancient Computers interesting and useful,

-Steve Stephenson

July 15, 2010

# Introduction

If you stare at an old mechanical calculator it just sits there. It does no computing and is, therefore, not a computer. When a person starts punching the keys and turning the crank the person-device computes and is a computer. So too, an abacus is just an assembly of beads and rods or lines and pebbles, and is not a computer. But when a person uses the abacus to perform calculations, then the person-abacus is a computer; a remarkably fast and accurate computer, as demonstrated by a Japanese abacus (Soroban) operator who beat a skilled electric calculator operator in a contest in Tokyo on Nov 12, 1946. (Kojima, 1954, p. 12)

In classic Greek architecture, an abacus is a flat slab of marble on top of a column's capital, supporting the architrave, or beam. Such an abacus (perhaps chipped beyond use in construction) makes a fine flat surface on which to inscribe lines; from which we get the name, counting board or line abacus[3]. Developed later, constrained bead devices with less arithmetic functionality[4] are also called abaci[5], e.g., Roman Hand Abacus, Chinese Suan Pan, and Japanese Soroban.

The Ancient Romans were excellent practical engineers and architects. Even today we marvel over their accomplishments and wonder how they did them. For example, in a BBC2 sponsored series Building the Impossible[6], Episode 2: The Roman Catapult, structural engineer Chris Wise wonders how the Romans did the calculations necessary to design and build the Roman Catapult used to destroy the walls of Jerusalem in 70 AD[7], when the math necessary wouldn't be developed for over 1500 years![8].

> “The Roman expression for 'to calculate' is 'calculos ponere' - literally, 'to place pebbles'. When a Roman wished to settle accounts with someone, he would use the expression 'vocare aliquem ad calculos' - 'to call them to the pebbles.'” (Jen, 1998; Menninger, 1969, p. 316)

Certainly the Romans would also use their abaci for engineering calculations[9]. Indeed, in his reference 16, Prof. Netz (2002) writes, "... my guess is that, mental calculating prodigies aside, complicated calculations were always done with the abacus."

Historians have published conjectures for what ancient counting board abaci looked like and how they were used (Ifrah, 2000, pp. 200-211; Lang, 1964; Menninger, 1969, pp. 295-306). On p.205, Ifrah concludes, "Calculating on the abacus with counters was ... a protracted and difficult procedure, and its practitioners required long and laborious training."

That is not true.

## Clues

Clues for the true structure and methods of the ancient counting board abaci are contained in three extant artifacts: The Japanese Soroban (Kojima; Menninger, pp.307-310), The Roman Hand Abacus (Ifrah, p.210; Menninger, p.305), and The Salamis Tablet (Ifrah, p.201; Menninger, p.299). Citing the book Mathematical Treatises by the Ancients, compiled by Hsu Yo toward the close of the Later Han dynasty (A.D. 25-220), Japanese historians of mathematics and the abacus corroborate the existence of the Roman Hand Abacus (See Appendix C: Chinese Abacus).

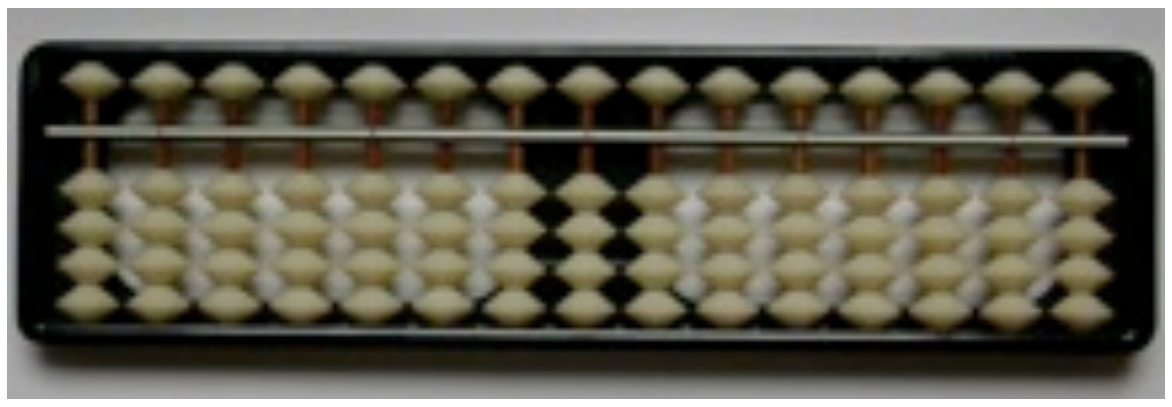

*Figure 1: A Small Japanese Soroban*

Every rod on a Soroban represents one decimal digit (Fig. 1). The bead above the bar represents five of the beads below the bar. Each rod can count from zero, no beads next to the bar, to nine, all beads next to the bar.

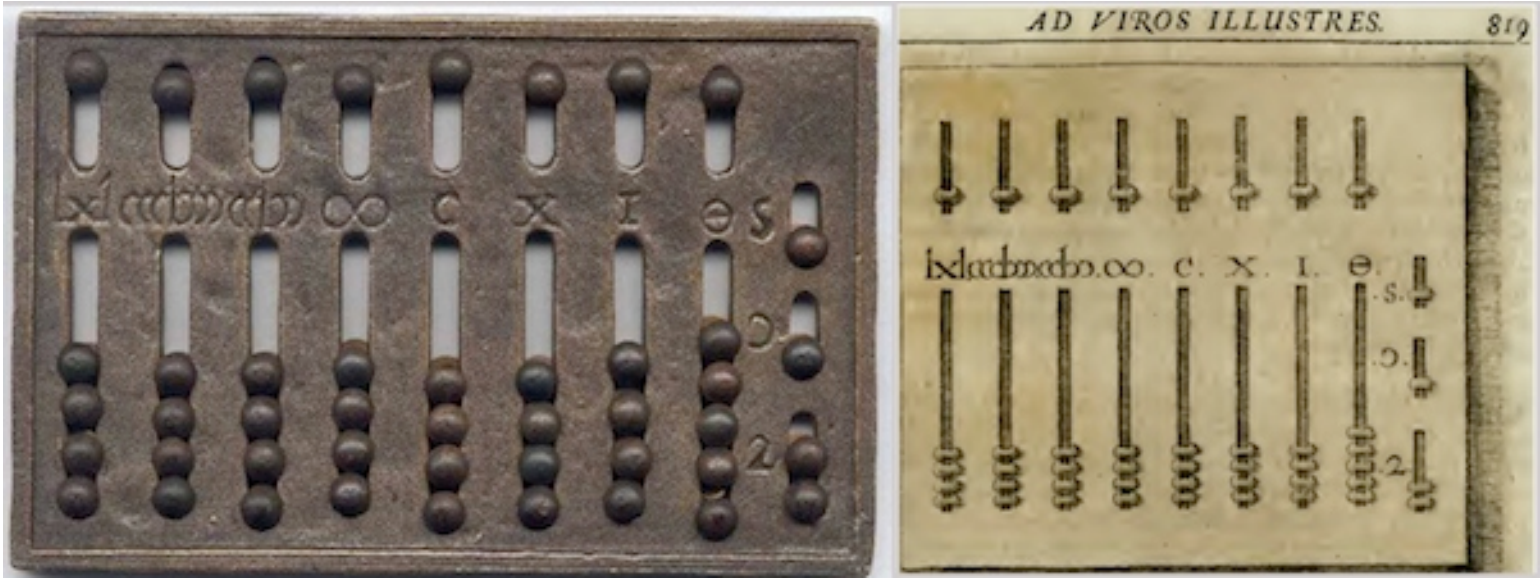

*Figure 2: Roman Hand Abacus*
*Left: Replica (Lutjens); Right: Source (Welser, p.422 & p.819)*

On the Roman Hand Abacus (Figure 2), each of the seven decimal digits has four beads in the lower slot and one bead in the upper slot; functioning exactly as the Soroban. It would be hard to understand why the Romans would not have developed similarly efficient methods to use the Hand Abacus as the Japanese did to use the Soroban (Kojima).

The two rightmost columns handle the Roman's base-12 fractions and both count to twelve, but differently. The left column counts to five in the lower slot and carries into the upper slot on a six count, repeats to a count of eleven, then carries into the decimal units column on a twelve count. But the rightmost column breaks each six count into two three counts. Why the difference?

Menninger, p.315, writes,

> "For more extensive and complicated calculations, such as those involved in Roman land surveys, there was, in addition to the hand abacus, a true reckoning board with unattached counters or pebbles. The Etruscan cameo and the Greek predecessors, such as the Salamis Tablet and the Darius Vase, give us a good idea of what it must have been like, although no actual specimens of the true Roman counting board are known to be extant."

Let us assume the obvious, that the Roman counting board abacus was The Salamis Tablet. Mapping the Roman Hand Abacus slot symbols onto The Salamis Tablet of Figure 3 results in Figure 4.

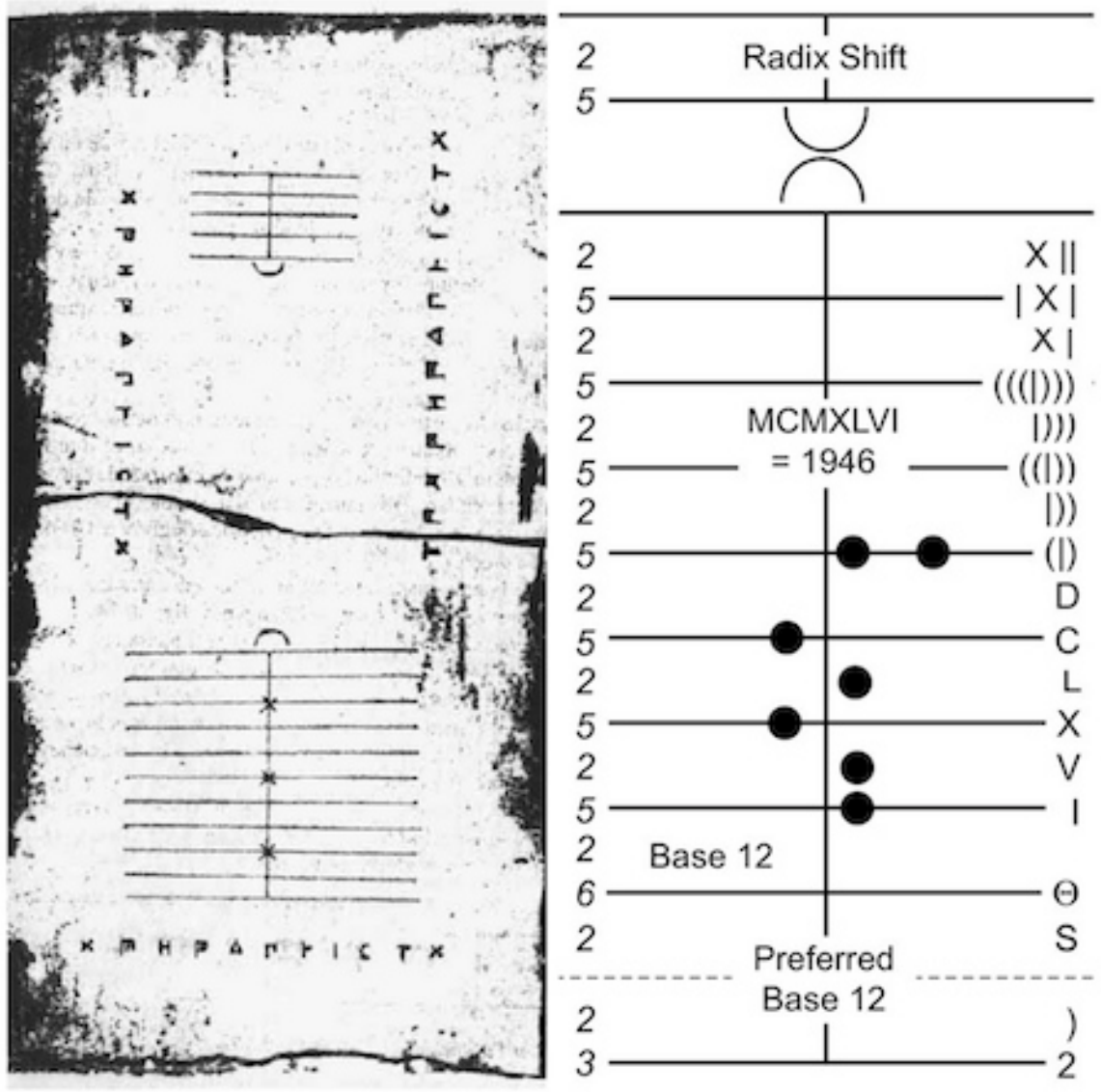

*Left: Figure 3, The Salamis Tablet (c. 300BC, Ifrah p.201)*
*Right: Figure 4: Roman Hand Abacus mapped onto Salamis Tablet*

The mapping is perfect. It uses the bottom grid's eleven lines exactly, no more, no less. However, to do so the Romans had to use a less preferred structure for one of the base-12 digits. Why it's less preferred is addressed below. But the fact that they had to make an engineering compromise is indicative that they used The Salamis Tablet as a design template for their Hand Abacus.

Note that both the bottom and top grids of The Salamis Tablet are divided into two equal left-and-right halves. This implies that there could be equal number of pebbles on both sides; a point ignored by historians when forming their conjectures about how The Salamis Tablet was used. (See also Appendix A: Arithmetic.)

In Figure 4, the numbers on the left are the promotion factors that are dictated by the Roman Numerals mapped to the lines and spaces. For example, a promotion factor of 5 means that five pebbles on that line can be replaced by one pebble in the space above. All the spaces between lines have a promotion factor of 2. (The unused dashed line is explained below.)

Placing the Roman Numeral MCMXLVI = 1946 on the mapped Salamis Tablet demonstrates the use of the left side of the grid as a subtractive side. This is an extremely important observation as it reduces the number of pebbles needed **tremendously**, both in total, and on each line or space, (see Appendix P: Pebble Efficiency). It makes the calculations much more efficient, robust, and rapid. An abacist would only need to carry a bag of pebbles that fits easily in a modern pants pocket to perform any of the four arithmetic calculations on any rational numbers with 10 decimal or 5 duodecimal (or sexagesimal) significant digits. (Stephenson, video 9.1)

As an example of using the subtractive side, if in the year 2009 = MMIX you wanted to calculate the age of a person born in 1946 = MCMXLVI, you would first make the 1946 in Figure 4 negative by moving each pebble to the opposite side of the vertical median line. (Follow along with pennies as pebbles and an abacus drawn on paper.[10]) Then to make room for the next addend, you would slide each of the pebbles as far away from the vertical median line as possible, left or right. Now you would add 2009 = MMIX by placing pebbles next to the median line: two on the right side of the "(|)", or M line, one to the left side of the "I" line, and one to the right side of the "X" line.

Merging the pebbles and replacing the C pebble with two L pebbles and one of the X pebbles with two V pebbles (both operations being demotions), then removing zero-sum pairs on every line and space (a pebble on each side of the median is a zero-sum pair) you read the answer LXIIV = 63.

While "IIV" is not considered by some people as a proper Roman Numeral expression, there are some rare examples of documents printed with both "IIV" and "IIX" type constructions (Handy, 2000). Proper form aside, these constructions do reduce the pebble count on counting board abaci, both in total and on each line and space. Using these constructions, any number can be registered with no more than two pebbles per line and one per space.

Moving the pebbles of an accumulated sum away from the median as far as possible, there will always be room near the median to place another additive or subtractive number. Before combining the pebbles the abacist can check the addend for accurate pebble placement without damaging the accumulated sum. This checking, or auditing feature adds greatly to the accuracy and robustness of these methods. After checking, the pebbles are combined and moved to the outside ready for the next number to be added. Constrained bead abaci like the Suan Pan and Soroban cannot do this addend checking.

## Design Compromise

The preferred base-12 digit configuration will never have more than two pebbles per line, but the other base-12 digit configuration will have three pebbles on a line for counts of 3 or 9. That is objectionable, both for pebble count per line or space as well as the psychological problem of handling two types of entities on the board. The more alike everything is, the faster and more accurately the operator can perform. The width of the counting board abacus would also have to be increased by 40% (14 vs. 10 pebble widths per line).

In article 26 of The Aqueducts of Rome, Frontinus writes,

> "... the inch ajutage, has a diameter of 1 1/3 digits. Its capacity is [slightly] more than 1 1/8 quinariae, i.e., 1 1/2 twelfths of a quinaria plus 3/288 plus 2/3 of 1/288 more."

In base-12 (using a semi-colon for the radix point and commas to separate digits), 1/8 is ;1,6 and (1+1/2)/12+3/288+(2/3)(1/288) is ;1,7,10. So Frontinus is saying 1;1,7,10 is slightly more than 1;1,6. Prof. Turner (2007) says that 1;1,7,10 results from the calculation of (1+1/3)^2 / (1+1/4)^2, where the squares are calculated before dividing.

The Romans would not have been able to do this calculation on the Roman Hand Abacus, nor on a Roman Hand Abacus mapped onto The Salamis Tablet. The number 1;1,7,10 has three base-12 fraction digits, and a fourth digit would need to be calculated for proper rounding. The Roman Hand Abacus only has two base-12 fraction digits.

So the Romans would have done the calculation on three coupled Salamis Tablets, each configured with 5 preferred-configuration base-12 digits, as in Figure 5. (But see Appendix U and Frontinus' Duodecimal Abacus.)

The pebbles in Figure 5 are placed for the first step in squaring 1 1/3. The radix shift (what we call an exponent of the base) is shown as a zero-sum pair for both the Multiplicand and Multiplier to indicate no shift, so the unit line is the top line of the bottom grid.

Frontinus' calculations, and others, are performed in a set of videos (Stephenson, videos 10.1-10.3, et al). It's much easier and more insightful to watch someone doing the calculations than reading about how to do them. However, there are simplified tabular examples of Multiplication and Division in Appendix M.

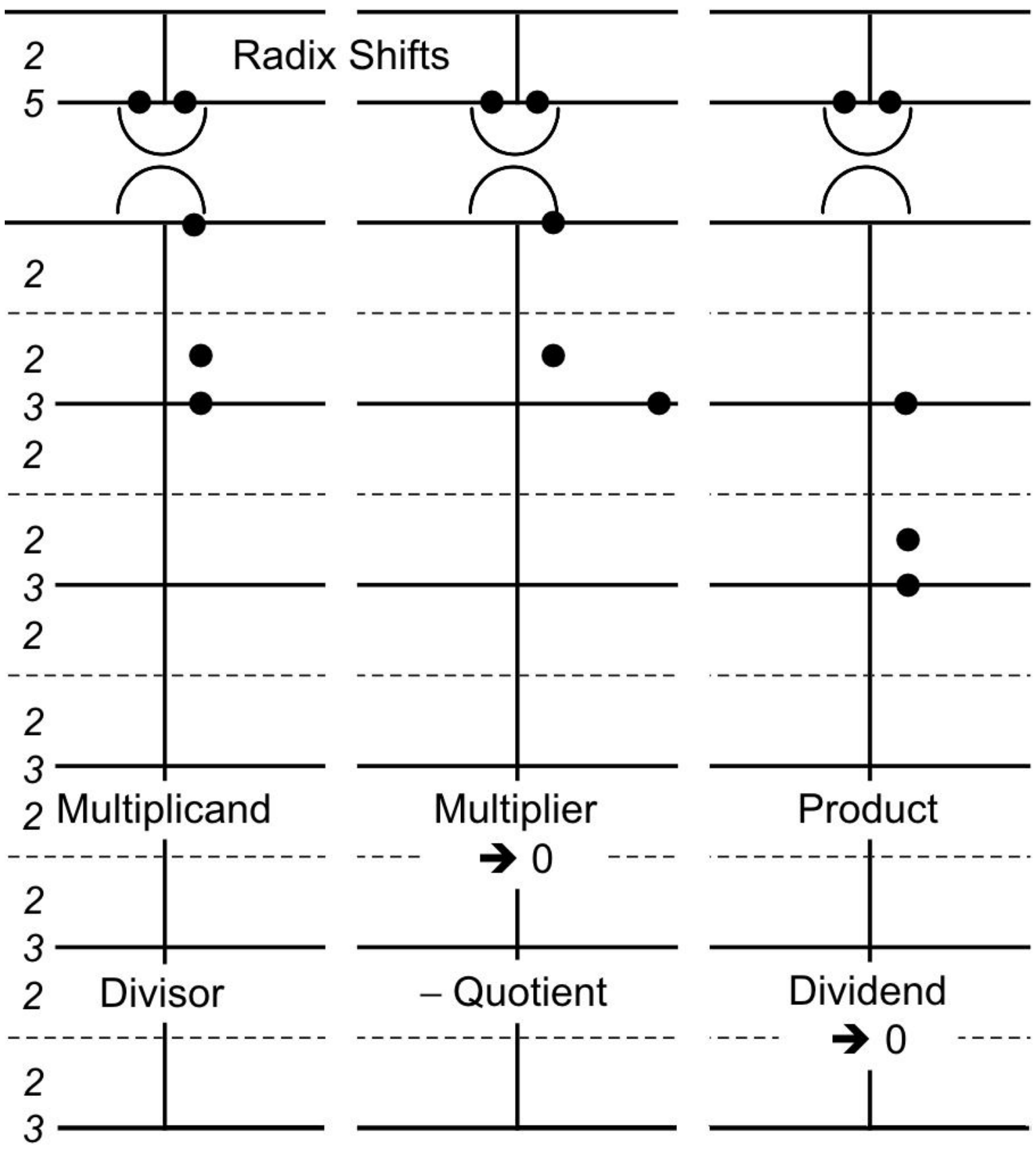

*Fig. 5. Three Coupled Base-12 Salamis Tablets for Multiplication or Division*

## Unused Dashed Lines

Why are the dashed lines unused? The Romans were borrowers. They borrowed The Salamis Tablet from the Greeks, but the Greeks borrowed it in turn from the Babylonians. The Babylonians used place value sexagesimal numbers written with reed styluses in cuneiform on clay tablets. Each of their digits contained 0-5 glyphs representing ten each and 0-9 glyphs representing one each, where there was at least one glyph; so each sexagesimal digit could represent the integers 1-59 (there was no zero symbol; see Figure 6).

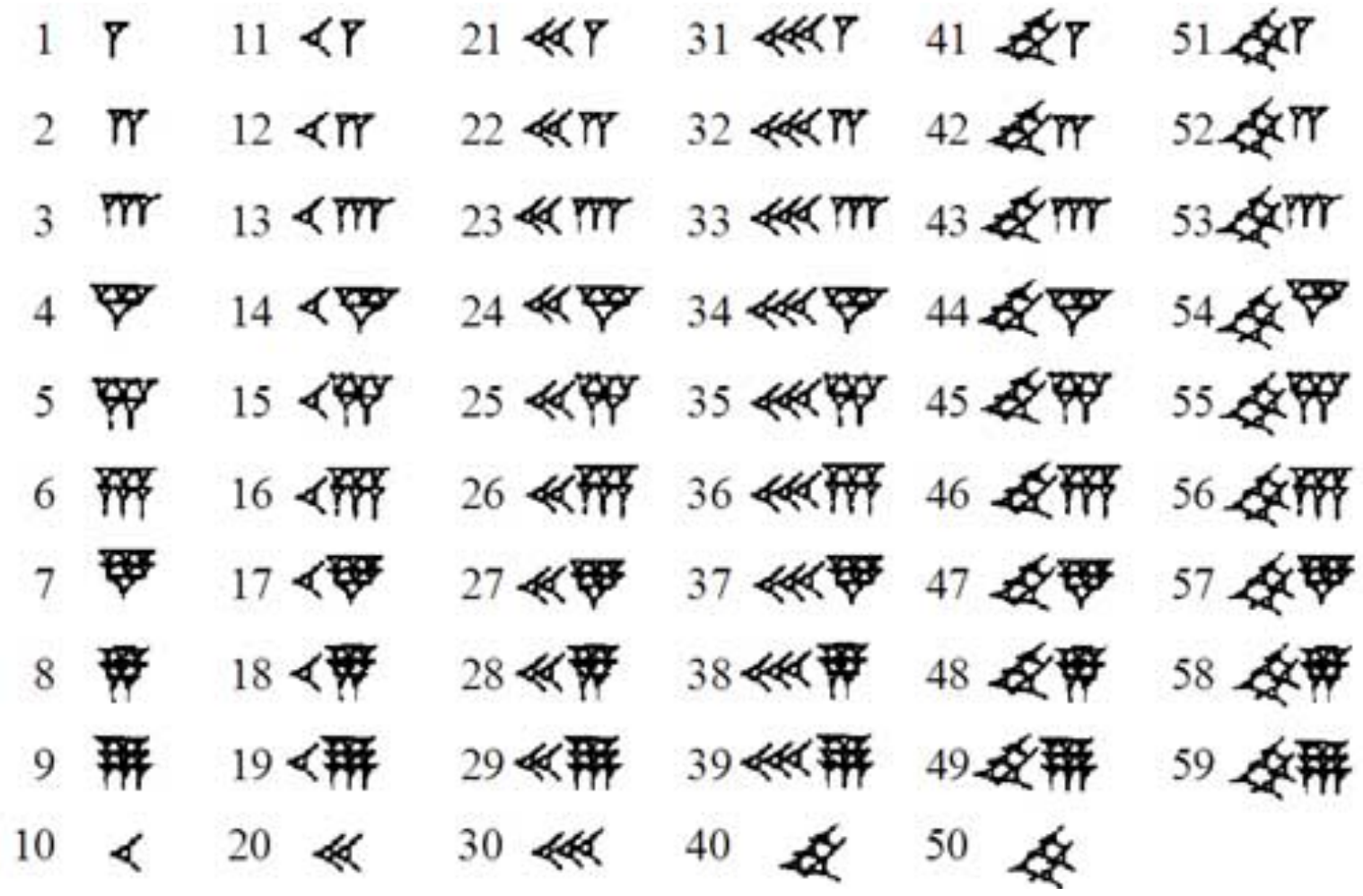

*Figure 6: Babylonian Numerals*

To use The Salamis Tablet with base-60 numbers every sequence of line, space, dashed-line, and space, from bottom to top are assigned promotion factors of 5-2-3-2. Notice that in both the Roman duodecimal and Babylonian sexagesimal configurations the value of the space below any line is one half the value above. In essence, the Romans used a subset of the sexagesimal system for their duodecimal system. (This argument is strengthened in Frontinus' Duodecimal Abacus.)

The Babylonians also did not have a radix symbol (decimal point). But the context was always specified so the radix shift was always known and the base-60 number was, therefore, always a fraction less than 1 and greater than or equal to 1/60. So why would you need a radix symbol (or trailing zeros, for that matter)?

Menninger pp.316-317 writes,

> " ... in our example above, the first step was to multiply the tens, 2 × 3; but in which column is the result, 6, to be placed? The rules of position in multiplication on counting boards have been puzzled over since time immemorial ... "

If all numbers are entered with the most significant digit in the top position, therefore as a fraction of one, along with the appropriate radix shift, then no complex positioning rules are needed. Just use the fact that

$$(M \times b^m) \times (N \times b^n) = (M \times N) \times b^{(m + n)}.$$

Archimedes is given credit for proving this exponent rule, but the Babylonians had probably been using it on their abaci for hundreds of years before.

The Babylonians lack of a radix symbol and elimination of complex positioning rules are strong evidence that the top grid on The Salamis Tablet is used for storage and manipulation of a radix shift, what we call an exponent of the base. Constrained bead abaci like the Suan Pan and Soroban do not have exponent capability.

So the Babylonians and perhaps the Sumerians had a computer using a place value number system and including the utility of exponential notation as early as 2300 BC[11]. Computers manipulating digits using digits of the hand, i.e., "digital computers" in use more than 4000 years ago!

But without a zero symbol, how could the Babylonians handle the one in twenty abacus results with empty digits (embedded zeros)[12]? Two ways come to mind:

1. Convert to another unit of measure; e.g., 104 yards = 312 feet; or

2. State the answer in multiple parts; e.g., 10403 yards = ten thousand, four hundred, and three yards. (Like writing a check.)

# Abaci Development

*a Plausible Historic Sequence*

## 1st abacus design

Your hands have 10 digits. 9834 with 24 pebbles:

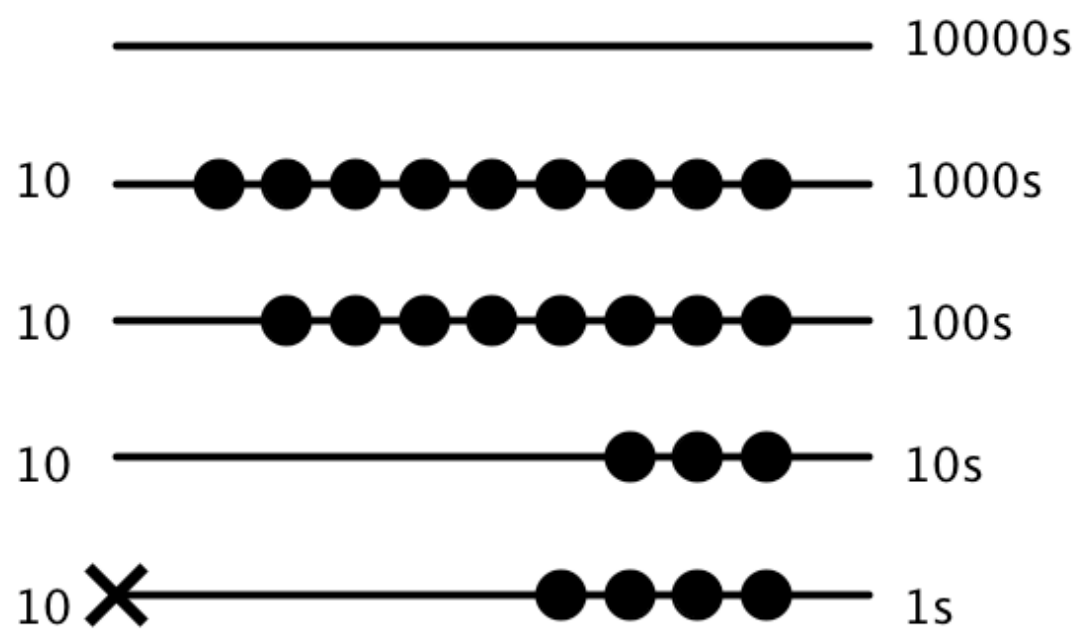

X marks the unit line; promotion factors along the left; pebble values along the right.

## 2nd abacus design

Each hand has 5 digits, 2 make ten. 9834 with 16 pebbles:

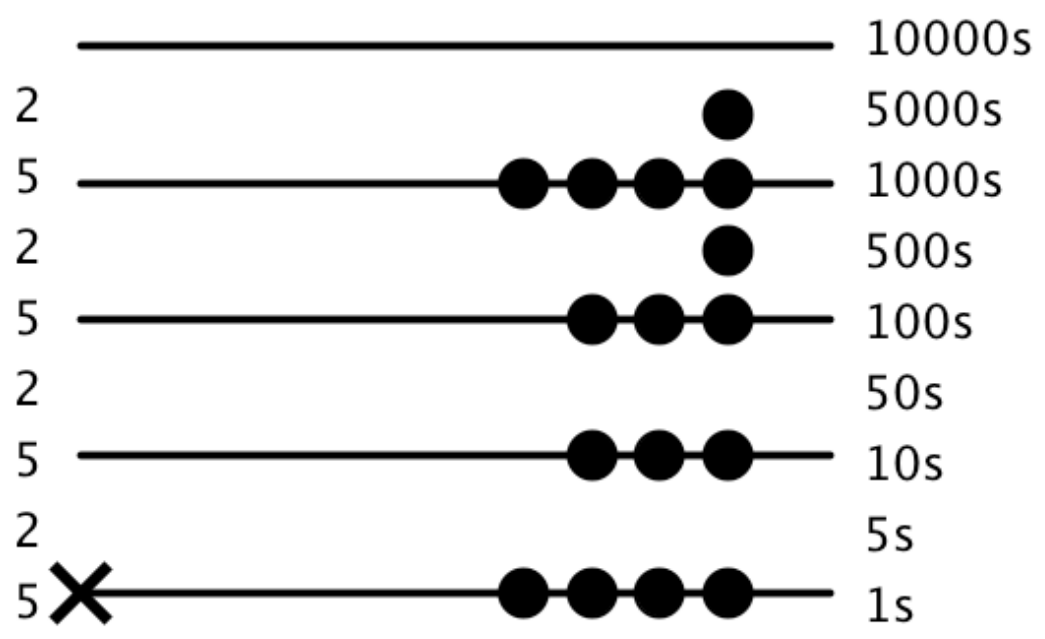

## 3rd abacus design

Opposites: yin-yang, male-female, left-right, etc. Since 9 = -1+10 = IX, 8 = -2+10 = IIX, 4 = -1+5 = IV, and 3 = -2+5; 9834 uses only 10 pebbles:

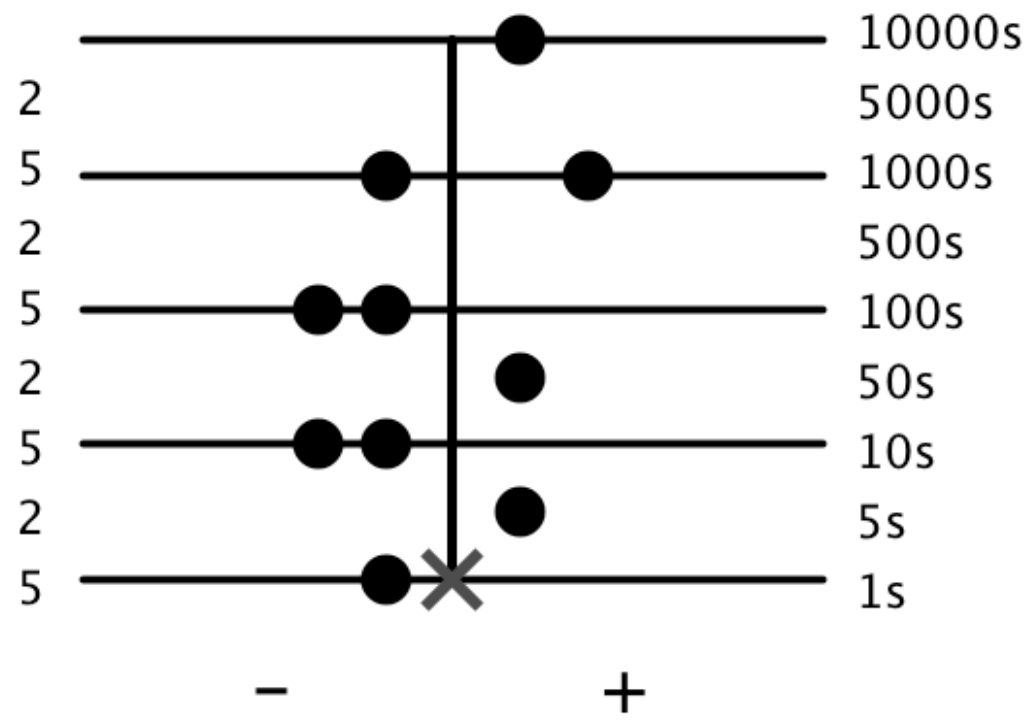

But -100 = (-50)+(-50) and -10 = (-5)+(-5), [demotions]:

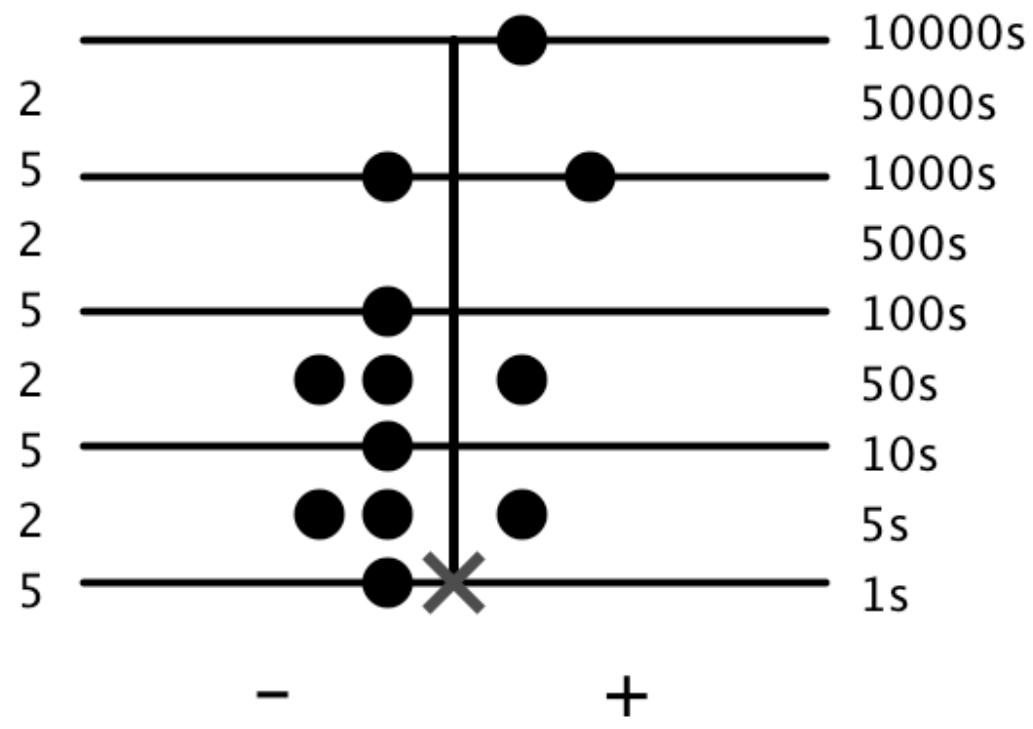

And since -k+k = 0 [cancellation of zero pairs]:

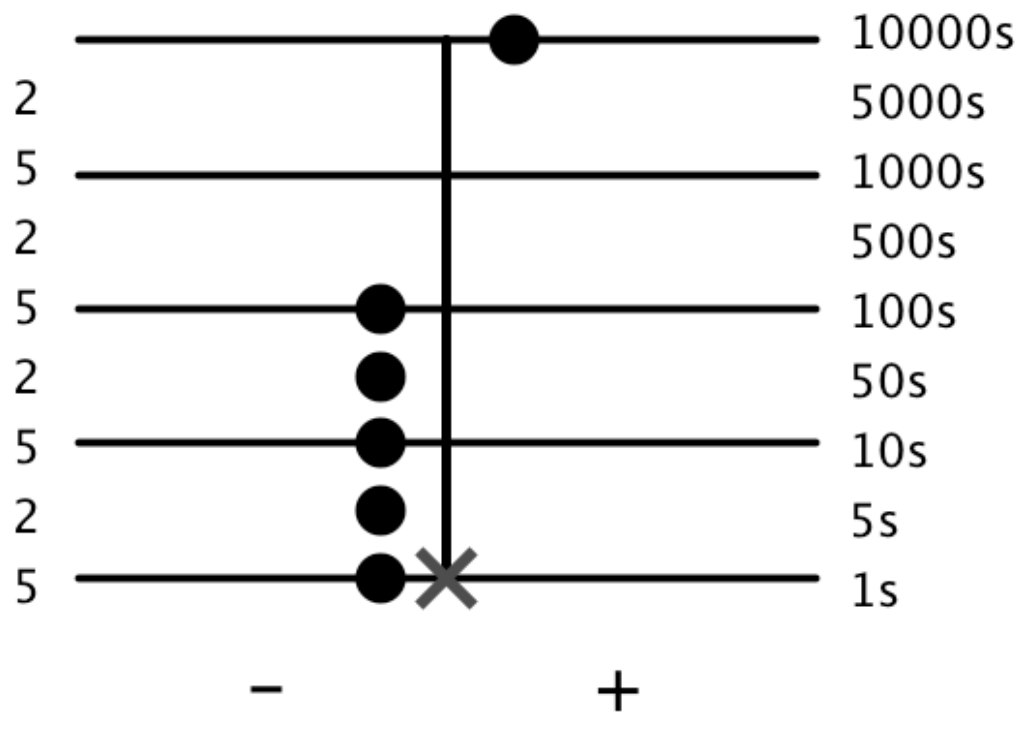

So 9834 can be represented with only 6 pebbles (9834 = - 166+10000 ).

## 4th abacus design

What if a number won't fit; e.g., 9,834,000,000,000,000 = 0.9834 x $10^{16}$?

To the ancients the exponent would be a scaling number[13]. They would enter 0.9834 at the top of the lower grid, then count how many lines down from the top the unit line would have been, and enter that count in a new upper grid. The unit line on the bottom grid is the top line. Only 9 pebbles are needed to represent 9,834,000,000,000,000:

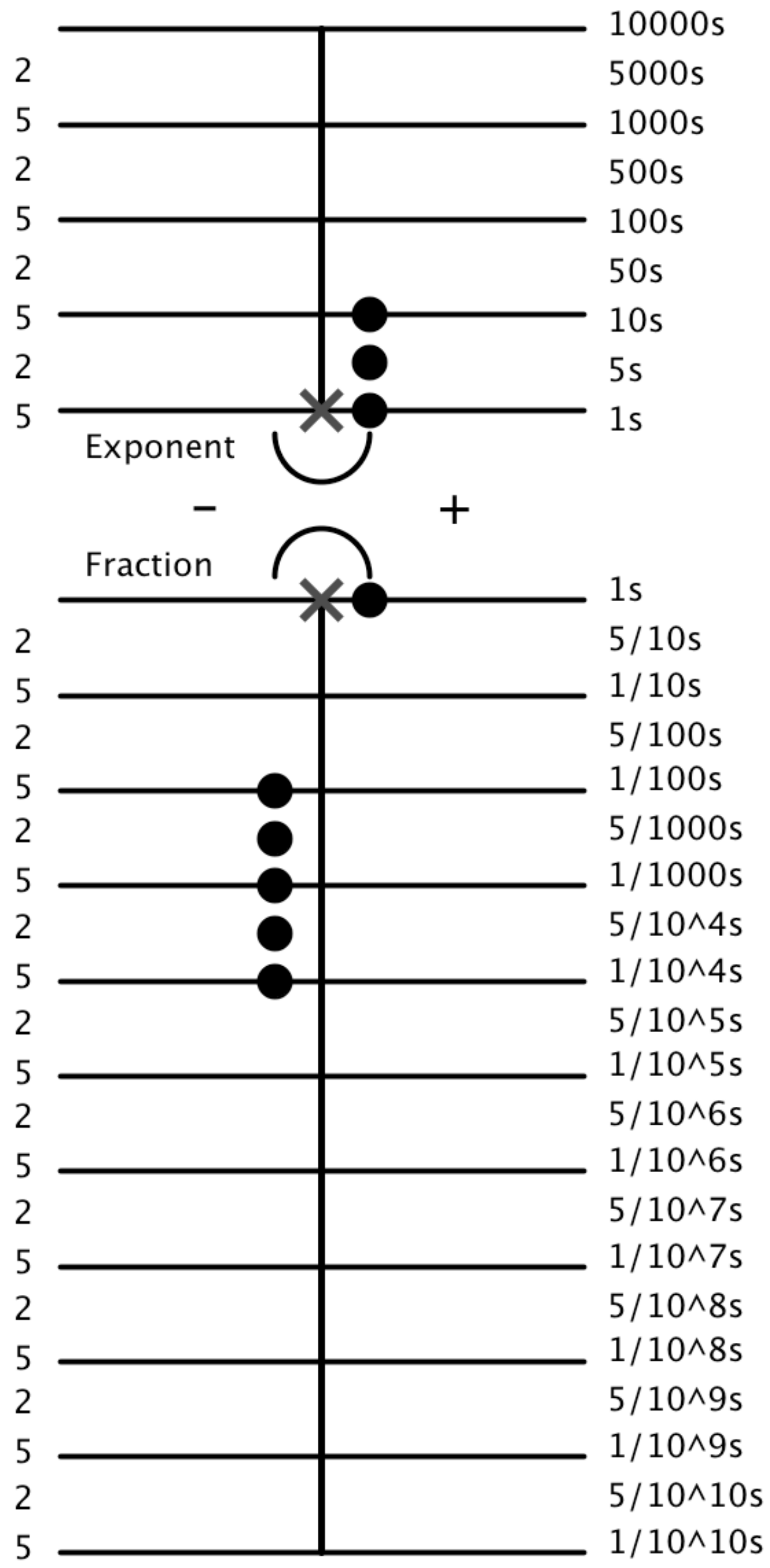

## 5th abacus design

60 = 10 × 6 = 5 × 2 × 3 × 2, so the abacus design is easily extended to sexagesimal, and the Babylonians' cuneiform numbers fit nicely into the structure; e.g., sqrt(2) = 1;24,51,10:

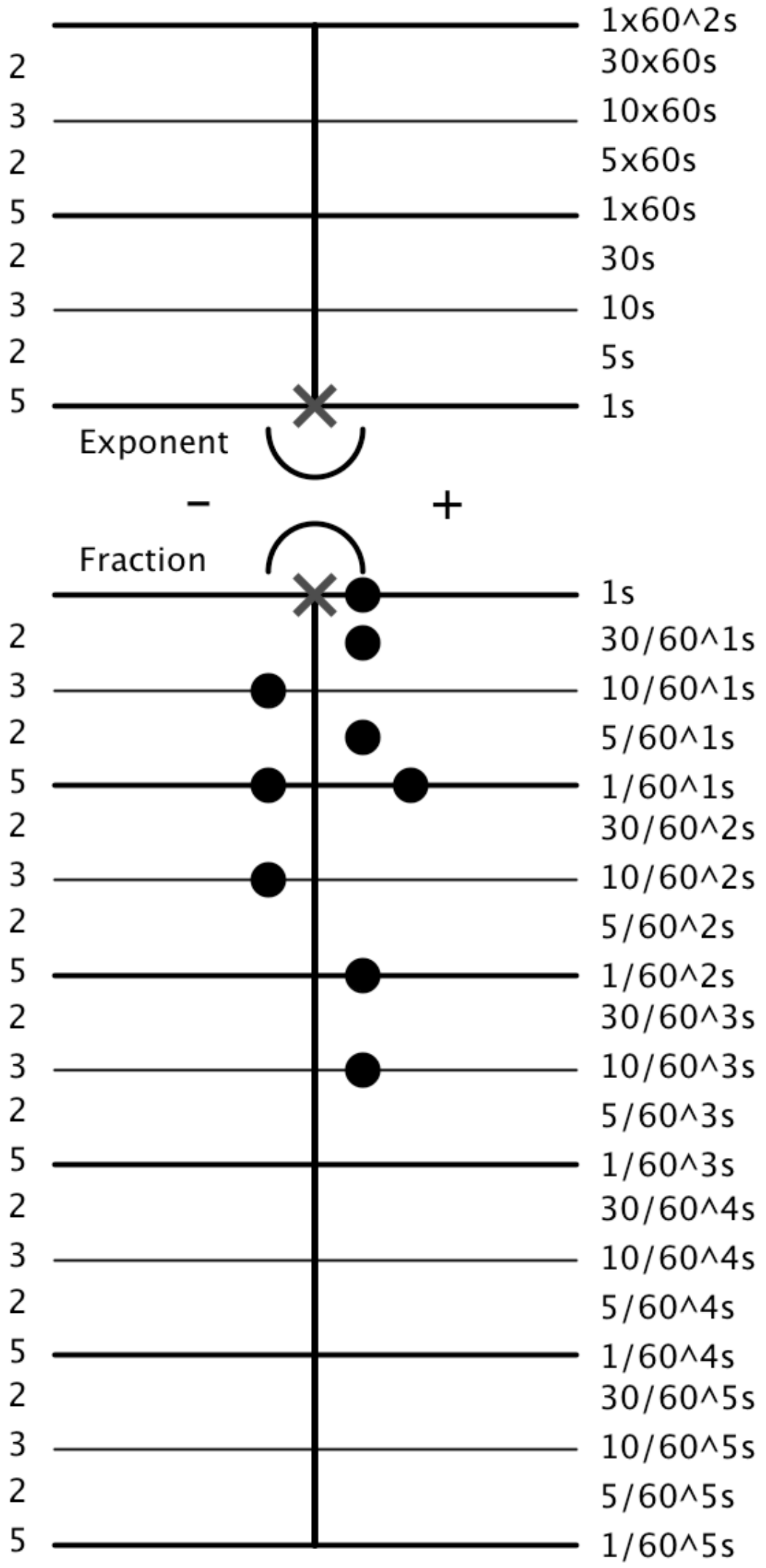

## Base-60

Why base-60, 360 degrees in a circle, 7 days in a week, etc.? All from Sumer around 2500 BC (Wilson, 2001), and inherited by the Babylonians. Wilson writes,

> "The Sumerians were great innovators in matters of time. It is to them, ultimately, that we owe not only the week but also the 60-minute hour. Such things came easily to people who based their maths not on a decimal system but on a sexagesimal one.
>
> Why were these clever chaps, who went for 60 because it is divisible by 2, 3, 4, 5, 6, 10, 12, 15, 20 and 30, fascinated by stubbornly indivisible seven? ...
>
> The Sumerians had a better reason for their septimalism. They worshiped seven gods whom they could see in the sky. Reverently, they named the days of their week for these seven heavenly bodies." (Wilson, 2001)

The divisibility of 60 was a convenient coincidental consequence, but not the primary reason the Sumerians adopted a sexagesimal number system[14]. They did so from the periods of the two slowest moving of their seven sky gods. Jupiter and Saturn take 12 and 30 years, respectively, to track through the Zodiac. The observant Sumerians knew this. The least common multiple of 12 and 30 is 60.

In 60 years Jupiter would go through 5 cycles and Saturn 2. We have 5 fingers on each of 2 hands. In both cases 5+2=7, the number of sky gods. The mystical Sumerians would think of this as manifestations of the sky gods reflecting themselves in our anatomy.

The product of 12 and 30 is 360, the number of degrees in a circle; did the Sumerians define the 360-degree circle? Probably, because dividing the Zodiac into 360 degrees means Jupiter traverses 30 degrees in a year and Saturn 12 degrees, thereby coupling the periods of the gods Jupiter and Saturn.

The Sun tracks through the Zodiac in one year. Jupiter would track 1/12 of the way in that time. Why not divide a year into twelfths, i.e., 12 months; then the Sun tracks the same distance in one month that Jupiter tracks in one year; thereby coupling the periods of Jupiter and the Sun. And since the Sun would then track 30 degrees along the Zodiac in a month, why not divide the month into about 30 days, the period of Saturn? Then the Sun tracks about 1 degree every day. Of course the Sumerians knew that a year is actually 365 days simply by watching the sun track through the Zodiac, so maybe they just added a 5 day Holiday, like the Egyptians did[15] … or, did the Egyptians copy the Sumerians?

## Salamis Tablet Designers

Both the astronomical and the anthropomorphic features of The Salamis Tablet in sexagesimal mode lead to the conclusion that the Babylonians, or their ancestors the Sumerians, were its designers; see Figure 7 (which is a frame in Stephenson's video 9.1). The Egyptians[16], Greeks, and Romans borrowed it for decimal and duodecimal calculations.

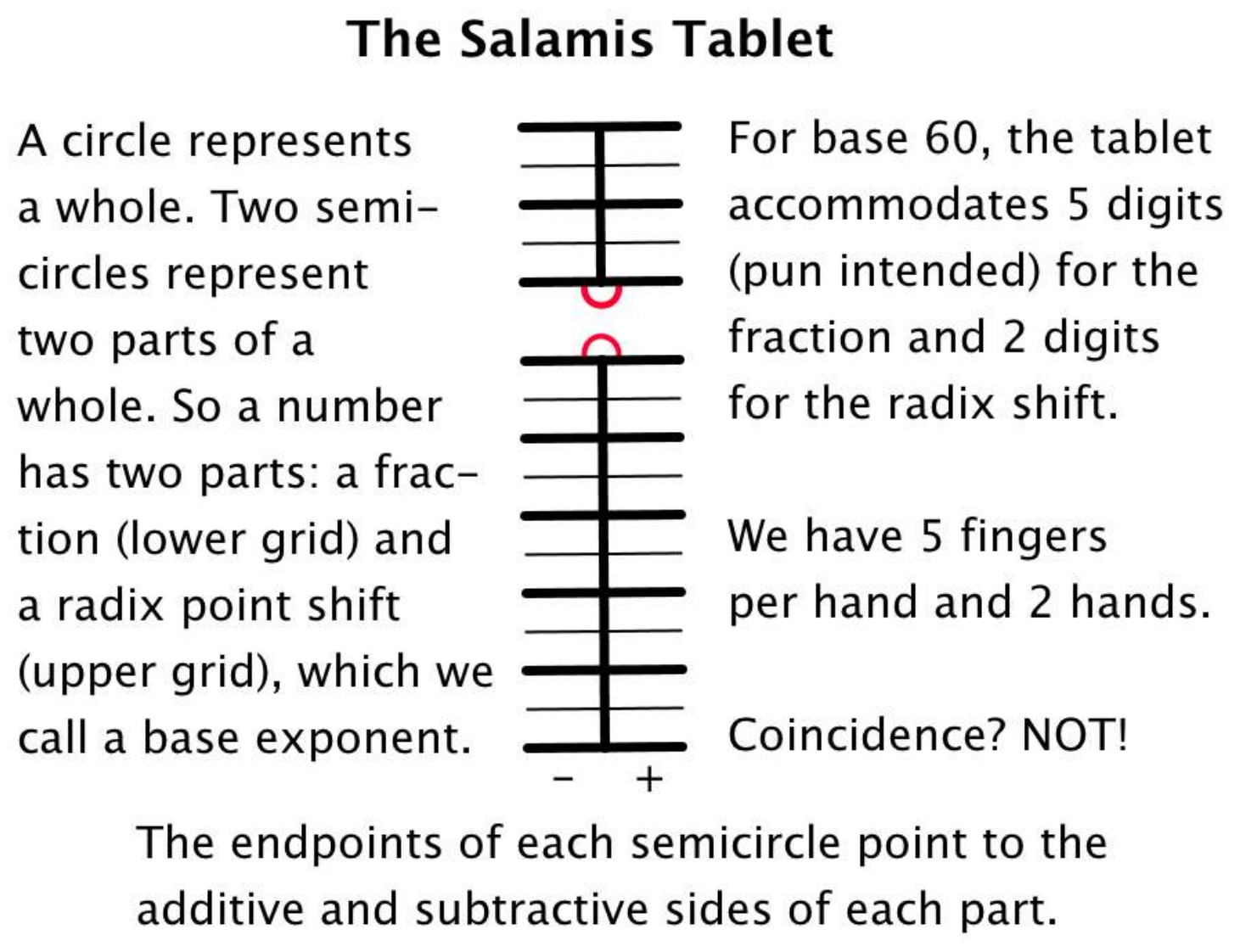

*Figure 7: Structural Design of The Salamis Table*

Robson (2008, p. 353, note 56) writes, "Both Proust 2000 and Hoyrup 2002a have argued, from numerical errors in Old Babylonian mathematical problems, that a place-marked abacus board with five columns [digits] was in use by the early second millennium [BC] if not earlier. This topic is in desperate need of careful and extended study..."

Perhaps this paper will aid that study, but certainly there are scholarly arguments strongly indicating that The Salamis Tablet was being used by the Old Babylonians by about 2000 BC or before.

## Computing in Sexagesimal

Some assyriologists, using modern pen, paper, and computers to produce beautiful multi-font typeset manuscripts, and using modern algebra and arithmetic as their basis for analysis, theorize that the Babylonians used tables of previously calculated values, recorded on clay tablets, to do their calculations.

For example, The MacTutor History of Mathematics web page states:

> "Perhaps the most amazing aspect of the Babylonian's calculating skills was their construction of tables to aid calculation. Two tablets ... give squares of the numbers up to 59 and cubes of the numbers up to 32. ... The Babylonians used the ... formula ab = [ (a + b)^2 - (a - b)^2 ] / 4 which shows that a table of squares is all that is necessary to multiply numbers, simply taking the difference of the two squares that were looked up in the table then taking a quarter of the answer." (O'Connor, 2008)

To use this formula to multiply two five digit sexagesimal numbers, proceed as follows:

(;a1,a2,a3,a4,a5)(;b1,b2,b3,b4,b5)
= (;a1)(;b1,b2,b3,b4,b5)
+ (;0,a2)(;b1,b2,b3,b4,b5)
+ (;0,0,a3)(;b1,b2,b3,b4,b5)
+ (;0,0,0,a4)(;b1,b2,b3,b4,b5)
+ (;0,0,0,0,a5)(;b1,b2,b3,b4,b5)

= (;a1)(;b1)+(;a1)(;0,b2)+(;a1)(;0,0,b3)+(;a1)(;0,0,0,b4)+(;a1)(;0,0,0,0,b5)
+(;0,a2)(;b1)+(;0,a2)(;0,b2)+(;0,a2)(;0,0,b3)+(;0,a2)(;0,0,0,b4)+(;0,a2)(;0,0,0,0,b5)
+(;0,0,a3)(;b1)+(;0,0,a3)(;0,b2)+(;0,0,a3)(;0,0,b3)+(;0,0,a3)(;0,0,0,b4)+(;0,0,a3)(;0,0,0,0,b5)
+(;0,0,0,a4)(;b1)+(;0,0,0,a4)(;0,b2)+(;0,0,0,a4)(;0,0,b3)+(;0,0,0,a4)(;0,0,0,b4)+(;0,0,0,a4)(;0,0,0,0,b5)
+(;0,0,0,0,a5)(;b1)+(;0,0,0,0,a5)(;0,b2)+(;0,0,0,0,a5)(;0,0,b3)+(;0,0,0,0,a5)(;0,0,0,b4)+(;0,0,0,0,a5)(;0,0,0,0,b5)

= (1;)[(;a1) (;b1)]
+(;1)[(;a1)(;b2)+(;a2)(;b1)]
+(;0,1)[(;a1)(;b3)+(;a2)(;b2)+(;a3)(;b1)]
+(;0,0,1)[(;a1)(;b4)+(;a2)(;b3)+(;a3)(;b2)+(;a4)(;b1)]
+(;0,0,0,1)[(;a1)(;b5)+(;a2)(;b4)+(;a3)(;b3)+(;a4)(;b2)+(;a5)(;b1)]
+(;0,0,0,0,1)[(;a2)(;b5)+(;a3)(;b4)+(;a4)(;b3)+(;a5)(;b2)]
+(;0,0,0,0,0,1)[(;a3)(;b5)+(;a4)(;b4)+(;a5)(;b3)]
+(;0,0,0,0,0,0,1)[(;a4)(;b5)+(;a5)(;b4)]
+(;0,0,0,0,0,0,0,1)[(;a5)(;b5)]

If only five significant figures are needed in the result, the last three lines can be discarded. Then there are nineteen products to calculate using the formula a(i)b(j) = [( a(i) + b(j) )^2 - ( a(i) - b(j) )^2]/4. The result of (a(i) + b(j))^2 will be of the form (1;c(i,j)) half the time, in which case (a(i) + b(j))^2 = [(1;) + (; c(i,j))]^2 = (1;) + 2(; c(i,j)) + (; c(i,j))^2.

On average, the calculation of each partial product would require: two additions, two table lookups, half a doubling, two subtractions, and two halvings. Combining the 19 partial products will require another 18 additions, being careful to add into the proper place value. In all, there are 56 additions, 38 table lookups, 8 doublings, 38 subtractions, and 38 halvings; a total of 178 operations!

How would you keep track of all this if you were limited to using reeds to write cuneiform on clay tablets? How many errors would you make? How would you find them?

Wouldn't it be simpler to suggest that the Babylonians developed and used abaci, with built-in error checking, to do their calculations?

# YBC 7289 sqrt(2) Calculation

"The Babylonian clay tablet YBC 7289 (c. 1800-1600 BCE) gives an approximation of sqrt(2) in four sexagesimal figures" as 1;24,51,10, which is approximately 1.414212963 in decimal, a remarkably accurate achievement for the time (0.0000423% too small).[17]

The calculation of sqrt(2) on a set of sexagesimal Salamis Tablets using Heron's Method (see Figure 8) takes 25 minutes in Stephenson's videos 8.1-8.3.

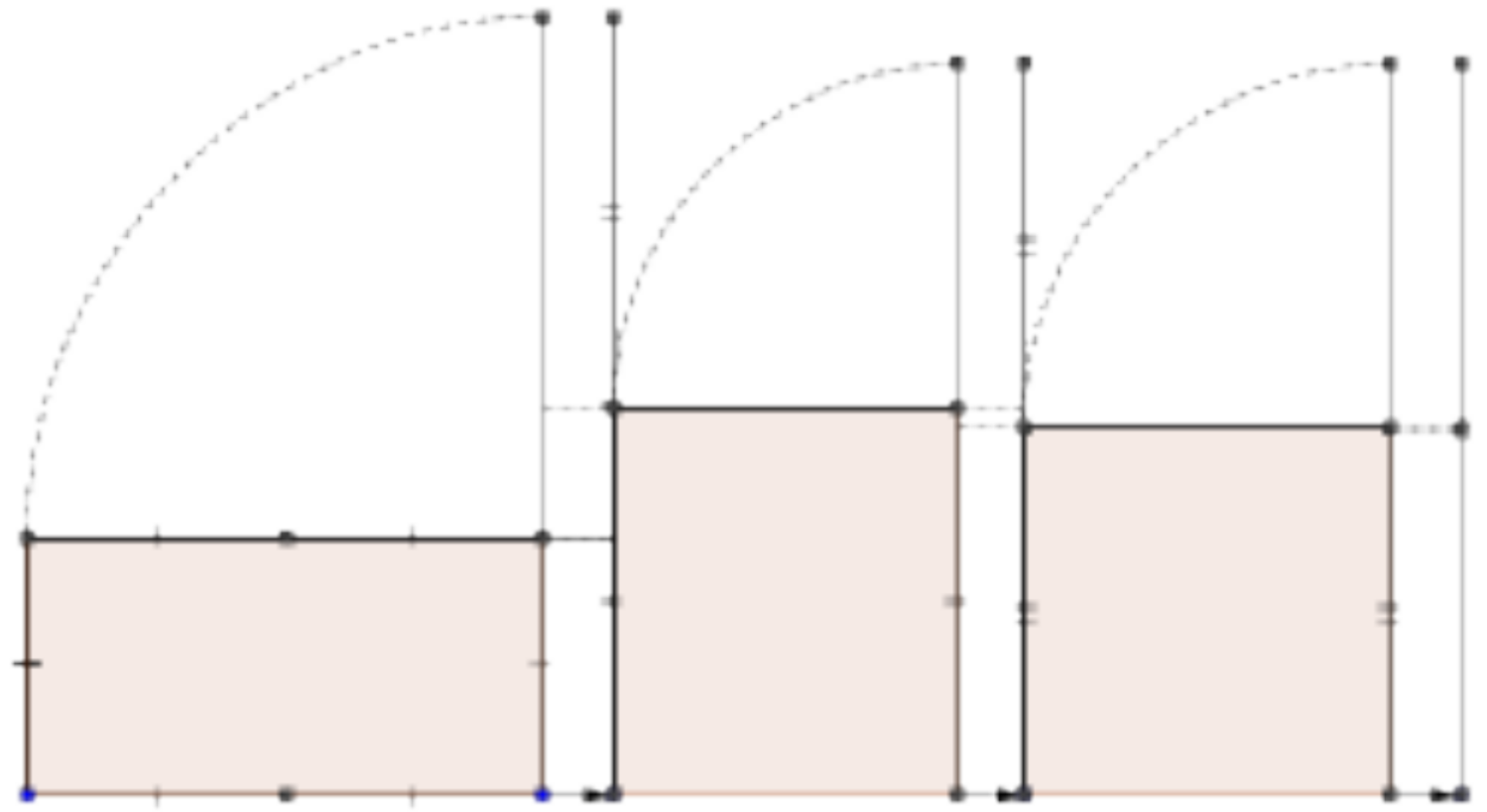

*Figure 8: Using Heron's Method to find sqrt(2)*
*The shaded rectangles have equal areas. See the long division symbols?*

Heron's method in algebra:

a(0) = 1, b(0) = 2;  a(1) = (a(0) + b(0)) / 2, b(1) = b(0) / a(1);  ...
a(i) = (a(i-1) + b(i -1)) / 2, b(i) = b(0) / a(i);

When b(i) ~= a(i), then b(i) ~= sqrt(b(0)) ~= sqrt(2).

If Historians performed the methods they think the Babylonians used (table lookups and such), and did so using only clay and reed styluses, just like the Babylonians, they would probably take longer than Stephenson's videos ... a LOT longer. The Historians should record their own performance video and post it to YouTube.com so we can compare its length to Stephenson's.

A possible analogy: in their early arithmetic practice exercises, modern elementary students learn how to perform long division and multiplication, and are required to show their work on paper. Then they are given four function electronic calculators and suddenly all of their practice exercises show only the problem statement and the answer, with no work.

Just like the mysterious school tablets of the Old Babylonian period![18]

## Conclusions

Features of the Roman Hand Abacus indicate that the Romans used a counting board abacus exactly like The Salamis Tablet for their heavy-duty calculations; it also gives us the promotion factors between lines and spaces. The Subtractive Notation of Roman Numerals indicates that one side of The Salamis Tablet grids are used for the additive part of a number and the other side for the subtractive part of the number.

The requirement of a four digit base-12 fractional answer to the calculation in article 26 of the book, The Aqueducts of Rome, indicates that the Romans used The Salamis Tablet in multiple configurations: decimal or duodecimal.

From these observations a consistent and powerful set of operational methods that the Romans must have used on The Salamis Tablet abacus has been identified.

The linkage between promotion factors of Roman duodecimal and Babylonian sexagesimal configurations, as well as astronomical and anthropomorphic features in the structure of the sexagesimal configuration, lead to the conclusion that the Babylonians, or their ancestors the Sumerians, were the designers of The Salamis Tablet and that they used much the same operational methods as the Romans, albeit with a different base. The Babylonians' lack of a radix point symbol suggests use of the upper grid of The Salamis Tablet to store and manipulate a radix point shift, or scaling factor; what we call an exponent of the base.

Two important conclusions are reached:

1. The ancients had FAST and POWERFUL computers, with BUILT-IN ERROR CHECKING, to use to power their empires; and

2. Now we know what those computers looked like and how they were used to perform all four arithmetic operations on decimal, duodecimal, and sexagesimal numbers; numbers that could be represented in an exponential notation, where both the fractional part and the exponential part could be positive or negative. The decimal numbers could have up to 10 significant digits in their fractional part and 4 significant digits in their exponent part. For duodecimal or sexagesimal numbers the fraction and exponent parts had up to 5 and 2 significant digits.

# Appendix A: Arithmetic

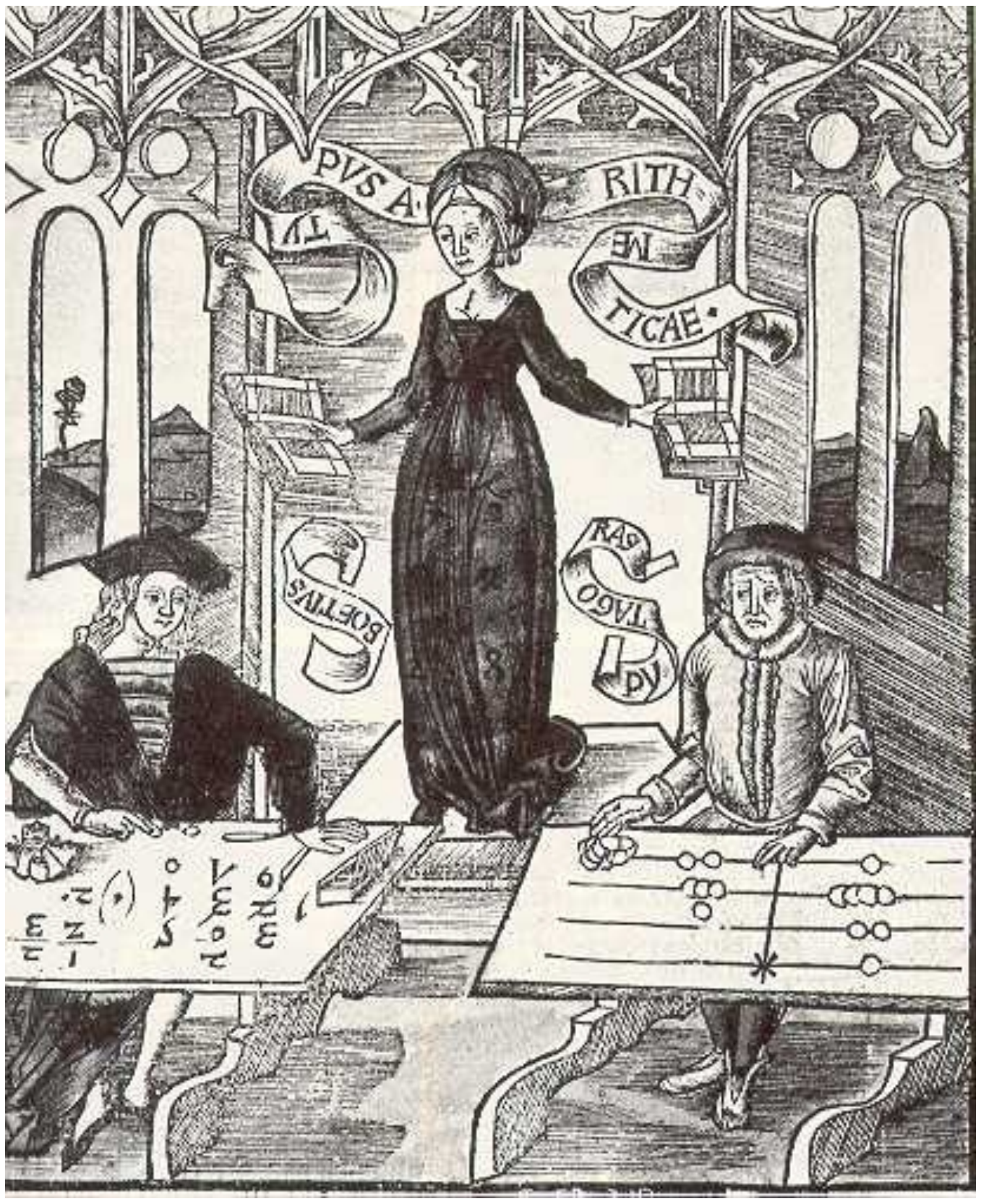

"This illustration of Arithmetic is from the *Margarita Philosophica* [written by Gregorius Reisch] which first came out in 1503. Arithmetic is illustrated by a beautiful lady holding a book in each hand. Sitting at two tables are Boethius and Pythagoras."[19]

Note that:

1. Pythagoras has his abacus oriented with a vertical median line;

2. The lines are equally divided so the same number of jettons can be accommodated on either side of the median line;

3. The top horizontal line is marked with an X to indicate the unit line; and

4. On the right side of the abacus is a jetton in a space between horizontal lines; all other jettons are on horizontal lines.

Prof. Barnard (1916, pp. 313-314) describes and quotes from Legendre, François, 1753, L'Arithmétique en sa perfection, Paris, pp. 497-528, Trait_ de l'arithmetiqué par les jetons:

> It was permissible to set and to work the jettons of the sum without using the spaces [between lines], ... But it was much more convenient to anyone who was expert at the practice, and less confusing to the eye, to reduce the number of jettons by using the spaces.

# Appendix C: Chinese Abacus

In the first chapter of his second book, Kojima (1963) writes,

> "The only reliable account of the origin of the Oriental abacus is in a book entitled Mathematical Treatises by the Ancients compiled by Hsu Yo toward the close of the Later Han dynasty (A.D. 25-220) at the beginning of the third century and annotated by Chen Luan in the sixth century. This book gives some information about various reckoning devices of those days and was one of the Ten Books on Mathematics (Suan-hwei-shi-chu) which were included among the textbooks to be read for government service examinations in China and Japan for many centuries.
>
> Chen Luan in his note gives the following description of the calculating device: "The abacus is divided into three sections. In the uppermost and lowest section, idle counters are kept. In the middle section designating the places of numbers, calculation is performed. Each column in the middle section may have five counters, one uppermost five-unit counter and four differently colored one-unit counters."
>
> The extent to which the counting board was used may be told by Hsu Yo's poetical description of the board. The verse, which is highly figurative and difficult to decipher, may read: 'It controls the four seasons, and coordinates the three orders, heaven, earth, and man.' This means that it was used in astronomical or calendar calculations, in geodetic surveys, and in calculations concerning human affairs.
>
> The reader will notice a close similarity between this original Oriental abacus and the Roman grooved abacus, except for the difference that counters were laid down in the former while they were moved along the grooves in the latter. Because of this and other evidence, many leading Japanese historians of mathematics and the abacus have advanced the theory that the above-mentioned prototype of the abacus was the result of the introduction into the East of the Roman grooved abacus.
>
> The following corroborative pieces of evidence in favor of this theory are cited in the latest works by Prof. Yoemon Yamazaki and Prof. Hisao Suzuki of Nihon University:
>
> (1) The original Chinese abacus has a striking resemblance in construction to the Roman grooved abacus, as is evident in the foregoing quotation from Hsu Yo's book, e.g., four one-unit counters and one five-unit counter in each column.
>
> (2) The method of operation of the ancient Chinese abacus was remarkably similar to the ancient Roman method. [e.g.,] In ancient China, multiplication and division were performed by the repetition of addition and subtraction

(3) Traces of reckoning by 5's may be found in the Chinese pictorial representation of reckoning-block calculation as in the Roman numerals, as: six: VI (5 + 1) seven: VII (5 + 2) eight: VIII (5 + 3) four: IV (- 1 + 5)

(4) Trade was carried on between China and Rome. Chinese historical documents written in the Han dynasty (206 B.C.-A.D. 220) furnish descriptions of two land routes, called silk roads, connecting the two great empires. ...

Until the introduction of Western mathematics, mathematicians in China and Japan utilized reckoning-block calculation, which had not only been developed to the point of performing basic arithmetic operations but was also used to solve quadratic, cubic, and even simultaneous equations. It is presumed that they did not think it worthwhile to concern themselves with the other reckoning devices, including the abacus, which was, in their eyes, an inferior calculator barely capable of performing multiplication and division by means of the primitive cumulative method of addition and subtraction. Probably another reason which alienated mathematicians from these reckoning devices was that these instruments gave only the result of calculation, and were incapable of showing either the process of calculation or the original problem.

In ancient times China was primarily a nomadic and agricultural country, and business in those days had little need of instruments of rapid calculation. Anyway a millennium after the Han dynasty there was no record of the abacus. During the dozen centuries beginning with its first mention in the Han dynasty until its development, this primitive calculator remained in the background.

However, with the gradual rise of commerce and industry, the need for rapid calculation grew. The modern, highly efficient abacus, which probably appeared late in the Sung dynasty (906 - 1279), came into common use in the fourteenth century. The great rise and prosperity of free commerce and industry during the Ming dynasty (1368 - 1636) are presumed to have promoted the use and development of the abacus. A number of books on mathematics brought out in those days give descriptions of the modern Chinese abacus and give accounts of the modern methods of abacus operation, including those of multiplication and division." (Kojima, 1963)

So the Chinese DID NOT invent the Abacus. They DID copy the Roman Hand Abacus after trading with the Romans over the Silk Road. Then they dismissed it, thinking it an inferior device, and ignored it for a thousand years, until the rise of business required faster arithmetic calculation speed than their numerals or homegrown devices were capable of.

The mention of the abacus by Hsu Yo in 220 AD corroborates the existence of the Roman Hand Abacus.

## Appendix E: Egyptian-Babylonian Link

The following block-quotes are excerpts from Friberg (2005), interspersed with a couple of my comments prefixed with "sks: "

> “[p.27] ... in OB [Old Babylonian] cuneiform texts special cuneiform signs are used for the Babylonian basic fractions 3' (= 1/3), 2' (= 1/2), 3" (= 2/3), and 6" (= 5/6). In a similar way, special notations are used in [Egyptian] hieratic mathematical papyri for the hieratic basic fractions 6' (= 1/6), 4' (= 1/4), 3' (= 1/3), 2' (= 1/2), and 3" (= 2/3). All other fractions, not counting fractions of measures, are written as "parts" (also called "unit fractions") with dots over the numbers.”

sks: 6" and 3" are strong indicators of the use of subtractive notation on a counting board abacus. Here's how they would be represented:

*6" = 5/6 = 1 - 10/60:*

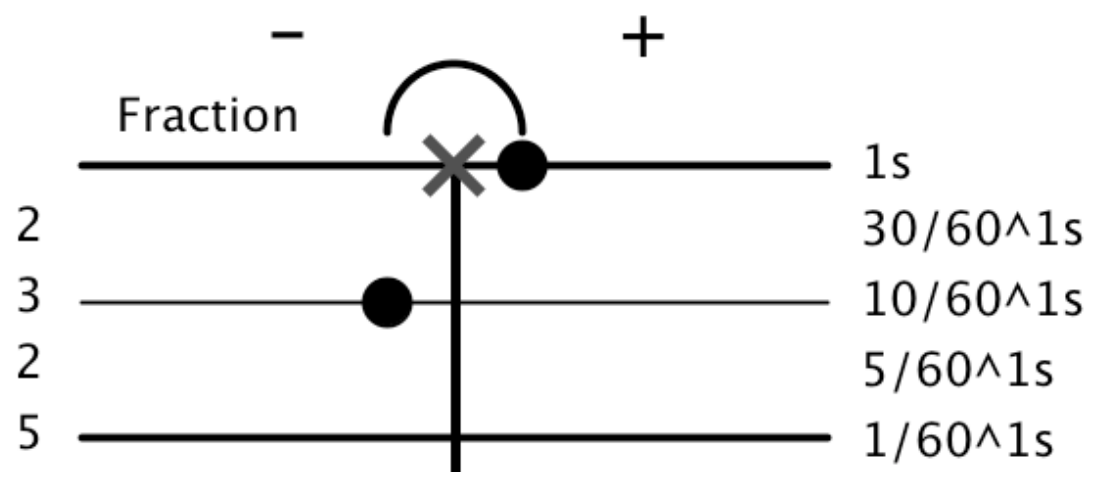

*3" = 2/3 = 1 - 20/60:*

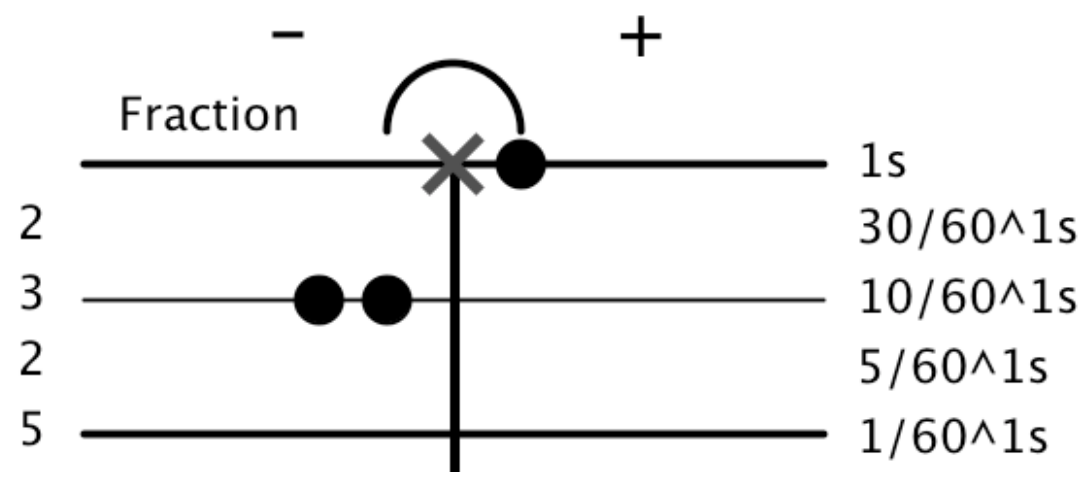

> “[p.32] Hence, if the quantity is 1 - 1/10 times 10 = 10 - 1 = 9. ... (Note that, apparently, in this exercise 9 divided by 10 is not represented by 3" + 1/5 + 1/30, a sum of parts, as in the 1/10 x n table in P.Rhind and in P.Rhind # 6, but by 1 - 1/10, a difference!)
>
> [p.113] In reference to the problem in Egyptian papyrus P.Cairo § 2 b (DMP # 5), "If it is said to you: 6" + 1/10 + 1/20 + 1/120 + 1/240 + 1/480 + 1/510, what remainder will complete 1?", Friberg writes, "The following tentative explanation ... is ... that while the author of P.Cairo, living in Hellenistic Egypt, nominally counted with traditional sums of parts [unit fractions] (traditional Egyptian), and with what looks very much like common fractions (a late Egyptian invention??), he may also have operated covertly with sexagesimal fractions (Babylonian)! ... [because] 6" + 1/10 + 1/20 +

1/120 + 1/240 + 1/480 can be equated to the sexagesimal fractions ;50 + ;06 + ;03 + ;00 03 + ;00 15 + ;00 07 30."

[p.128] In P.Cairo § 11 a (DMP # 36), ... Apparently, all calculations are carried out by use of sexagesimal arithmetic, although the results of the computations are expressed in terms of sums of parts.

[p.137] Referencing "m. P.Cairo §§ 15-16 (## "32-33")": "(Without the use of covert counting with sexagesimal numbers, the answer w = 8 3' 1/10 1/60 would have been instead that w = 8 3' 1/14 1/21.)"

[pp.165-166] Referencing P.British Museum 10794: "This small fragment contains only the first ten lines of two multiplication tables, one for 1/90, the other for 1/150. Note that 90 and 150 both can be expressed as regular sexagesimal numbers, 90 = 1 30 and 150 = 2 30. Their reciprocals are the sexagesimal fractions ;00 40 and ;00 24. The computation of the two tables probably made use of sexagesimal arithmetic. In his paper about Greek and Egyptian techniques of counting with fractions, Knorr, HM 9 (1982), 156, confessed that he (like Parker before him) was puzzled by 'idiosyncrasies' in the computational procedure ... Why was 2 x 1/150 given as 1/90 + 1/450 and not as 1/75, and why was 3 x 1/150 given as 1/60 + 1/300 and not as 1/50, and so on? ...The assumption that the scribe used sexagesimal arithmetic leads to a much simpler explanation ... . The ... slightly puzzling feature is why, in line 4, the sexagesimal fraction ;01 36 was not simply resolved as ;01 + ;00 36 = 1/60 + 1/100, and why similarly, in line 9, the sexagesimal fraction ;03 36 was not resolved as ;03 + ;00 36 = 1/20 + 1/100."

[pp.189-190] The rather detailed analysis ... above of the contents of P.Cairo clearly shows that in the third century BCE, if not sooner, Egyptian mathematics had become deeply influenced by Babylonian mathematics. ... Above all, the Babylonian influence is evident in the hidden use of sexagesimal fractions as a convenient computational tool, side by side with the Egyptian traditional use of sums of parts ...

[p.192] ... there can be little doubt that there were no significant differences between the general level and extent of the knowledge of mathematics in Egyptian demotic mathematical texts and in Mesopotamian cuneiform mathematical texts towards the end of the first millennium BCE, and that there are no signs of influence on either from high-level Greek mathematics.

[pp.269-270] ... the initial development of mathematical ideas started at a very early date with the invention of words for sexagesimal or decimal numbers in various ancient languages, and with the widespread use of number tokens in the Middle East. A major step forward was then the invention of an integrated family of number and measure systems, in connection with the inventions of writing in Mesopotamia and neighboring areas of Iran in the late fourth millennium BCE. There must have been a similar development in Egypt, about which not much is known at present. A small number of known examples of proto-Sumerian, Old Sumerian, Old Akkadian, and Eblaite mathematical exercises and table texts are witnesses of the continuing important role played by education in mathematics in the scribe schools of Mesopotamia throughout the third millennium [BCE].

Then there is a strange gap in the documentation, with almost no mathematical texts known from the Ur III period in Mesopotamia towards the end of the third millennium BCE. Nevertheless, at some time in the Ur III period a new major step in the development of mathematics was taken with the invention of sexagesimal place value notation. To a large part as a result of that invention, mathematics flourished in the Old Babylonian scribe schools in Mesopotamia. Simultaneously, mathematics may have reached a comparable level in Egypt, and, in spite of the fundamentally different ways of counting in the two regions, there was clearly some communication of mathematical ideas between Egypt and Mesopotamia.

A few late Kassite mathematical texts seem to indicate that the Old Babylonian mathematical tradition was still operative, although reduced to a small trickle, in the second half of the second millennium BCE.

Then follows a new strange gap in the documentation."

sks: Perhaps because of the eruption of Thera that destroyed the Minoan civilization, and much of the normal life in the greater eastern Mediterranean area.

"When mathematics flourished again in Mesopotamia in the Late Babylonian and Seleucid periods in the second half of the first millennium BCE, possibly in connection with the rise of mathematical astronomy, a great part of the Old Babylonian corpus of mathematical knowledge had been taken over relatively intact. However, for some reason, the transmission of knowledge cannot have been direct, which is shown by an almost complete transformation of the mathematical vocabulary.

Similarly in Egypt, after a comparable gap in the documentation, there was a new flourishing of mathematics, documented by demotic and Greek mathematical papyri and ostraca from the Ptolemaic and Roman periods. Some of the Greek mathematical texts are associated with the Euclidean type of high-level mathematics. Except for those, the remainder of the demotic and Greek mathematical texts show clear signs of having been influenced both by Egyptian traditions, principally the counting with sum of parts, and by Babylonian traditions. An interesting new development was the experimentation with new kinds of representations of fractions, first sexagesimally adapted sums of parts, soon to be abandoned in favor of binomial fractions, the predecessors of our common fractions.

The observation that Greek ostraca and papyri with Euclidean style mathematics existed side by side with demotic and Greek papyri with Babylonian style mathematics is important for the reason that this surprising circumstance is an indication that when the Greeks themselves claimed that they got their mathematical inspiration from Egypt, they can really have meant that they got their mathematical inspiration from Egyptian texts with mathematics of the Babylonian type."

# Appendix M: Multiplication and Division Examples

In Multiplication a partial product for every Multiplier pebble on a solid line is added into the Product at the proper place value, and the Multiplier pebble is taken away. If there are no more Multiplier pebbles on solid lines (representing powers of the base) then a Multiplier pebble in a space is replaced by the equivalent number of pebbles on the line below, or the Multiplier is doubled and the Multiplicand halved, or vice versa. When there are no more Multiplier pebbles the product is complete. The radix shift of the Product is the sum of the Multiplicand and Multiplier radix shifts.

For example, to multiply the decimal numbers 0.84 and 0.93:

| | Left Tablet | Middle Tablet | Right Tablet<br>Next to Median | <br>Away from Median |
|---|---|---|---|---|
| **Step** | **Multiplicand** | **Remaining Multiplier** | **Incremental Partial Product** | **Accumulated Partial Product** |
| 0 | 0.84 | 0.93 | 0 | 0 |
| 1 | 0.84 | 1.00 | 0.84 | 0 |
| | | -0.10 | | |
| | | 0.05 | | |
| | | -0.01 | | |
| | | -0.01 | | |
| 2 | 0.84 | -0.10 | -0.084 | 0.84 |
| | | 0.05 | | |
| | | -0.01 | | |
| | | -0.01 | | |
| 3 | 0.84 | 0.05 | | 0.756 |
| | | -0.01 | -0.0084 | |
| | | -0.01 | -0.0084 | |
| 4 | 0.42 | 0.10 | 0.042 | 0.7392 |
| 5 | 0.42 | | | 0.7812 |

In steps 1-4 the numbers in the Multiplier column are the values of individual pebbles. In step 4, instead of changing 0.05 to five 0.01s, the multiplier was doubled and the multiplicand was halved. Doubling and halving was practiced and documented extensively by the Egyptians.

In Division appropriate additive or subtractive copies of partial products of the Divisor are combined with the Dividend until the remaining Dividend is zero, accumulating the partial product pebbles in the middle grid as the Negative of the Quotient. When the Dividend is zero, simply move all the pebbles in the middle grid to the opposite side of the median and you have the Quotient. The Quotient's radix shift will be the radix shift of the Dividend added to the opposite of the radix shift of the Divisor.

For example, to divide the decimal number 0.84 by 0.93:

| | Left Tablet | Middle Tablet | | Right Tablet | |
|---|---|---|---|---|---|
| | | Next to Median | Away from Median | Next to Median | Away from Median |
| Step | Divisor | Accumulated Negative Quotient | Incremental Negative Quotient | Dividend Increment | Remaining Dividend |
| 0 | 0.93 | | | | 0.84 |
| 1 | 0.93 | | -1.0 | -0.93 | 0.84 |
| 2 | 0.93 | -1.0 | 0.1 | 0.093 | -0.09 |
| 3 | 0.93 | -0.9 | $-1/10^3$ | $-0.93/10^3$ | 0.003 |
| | | | $-1/10^3$ | $-0.93/10^3$ | |
| | | | $-1/10^3$ | $-0.93/10^3$ | |
| 4 | 0.93 | -0.903 | $-1/10^4$ | $-0.93/10^4$ | $0.21/10^3$ |
| | | | $-1/10^4$ | $-0.93/10^4$ | |
| 5 | 0.93 | -0.9032 | $-1/10^5$ | $-0.93/10^5$ | $0.24/10^4$ |
| | | | $-1/10^5$ | $-0.93/10^5$ | |
| | | | $-1/10^5$ | $-0.93/10^5$ | |
| 6 | 0.93 | -0.90323 | $1/10^6$ | $0.93/10^6$ | $-0.39/10^5$ |
| | | | $1/10^6$ | $0.93/10^6$ | |
| | | | $1/10^6$ | $0.93/10^6$ | |
| | | | $1/10^6$ | $0.93/10^6$ | |
| 7 | 0.93 | -0.903226 | $1/10^7$ | $0.93/10^7$ | $-0.18/10^6$ |
| | | | $1/10^7$ | $0.93/10^7$ | |
| 8 | 0.93 | -0.9032258 | $-1/10^8$ | $-0.93/10^8$ | $0.6/10^8$ |
| 9 | 0.93 | -0.90322581 | $1/10^9$ | $0.93/10^9$ | $-0.33/10^8$ |
| | | | $1/10^9$ | $0.93/10^9$ | |
| | | | $1/10^9$ | $0.93/10^9$ | |
| | | | $1/10^9$ | $0.93/10^9$ | |
| 10 | 0.93 | -0.903225806 | $-1/10^{10}$ | $-0.93/10^{10}$ | $0.42/10^9$ |
| | | | $-1/10^{10}$ | $-0.93/10^{10}$ | |
| | | | $-1/10^{10}$ | $-0.93/10^{10}$ | |
| | | | $-1/10^{10}$ | $-0.93/10^{10}$ | |
| 11 | 0.93 | -0.9032258064 | $-1/10^{11}$ | $-1/10^{11}$ | $0.5/10^{10}$ |
| | | | $-1/10^{11}$ | $-1/10^{11}$ | |
| | | | $-1/10^{11}$ | $-1/10^{11}$ | |
| | | | $-1/10^{11}$ | $-1/10^{11}$ | |
| | | | $-1/10^{11}$ | $-1/10^{11}$ | |
| 11 | 0.93 | -0.90322580645 | | | 0 |

So 0.84 / 0.93 = 0.9032258065. In the table, every number in the Incremental Negative Quotient column represents one pebble.

For both multiplication and division there is no need for reference tables, printed or memorized, and the calculations are blindingly fast compared to any other methods available to the Ancients.

# Appendix P: Pebble Efficiency

Assume two possibilities: using subtractive notation or not. If you have k pebbles, how high can you count on the 3rd abacus design? What patterns do the answers contain?

| | Without Negative Parts | | With Negative Parts | |
|---|---|---|---|---|
| k | Maximum Sequential Count | Increment | Maximum Sequential Count | Increment |
| 0 | 0 | 0 | 0 | 0 |
| 1 | 1 | 1 | 1 | 1 |
| 2 | 2 | 1 | 2 | 1 |
| 3 | 3 | 1 | 12 | 10 |
| 4 | 8 | 5 | 22 | 10 |
| 5 | 18 | 10 | 72 | 50 |
| 6 | 28 | 10 | 172 | 100 |
| 7 | 38 | 10 | 272 | 100 |
| 8 | 48 | 10 | 772 | 500 |
| 9 | 98 | 50 | 1,772 | 1,000 |
| 10 | 198 | 100 | 2,772 | 1,000 |
| 11 | 298 | 100 | 7,772 | 5,000 |
| 12 | 398 | 100 | 17,772 | 10,000 |
| 13 | 498 | 100 | 27,772 | 10,000 |
| 14 | 998 | 500 | 77,772 | 50,000 |
| 15 | 1,998 | 1,000 | 177,772 | 100,000 |
| 16 | 2,998 | 1,000 | 277,772 | 100,000 |
| 17 | 3,998 | 1,000 | 777,772 | 500,000 |
| 18 | 4,998 | 1,000 | 1,777,772 | 1,000,000 |
| 19 | 9,998 | 5,000 | 2,777,772 | 1,000,000 |
| 20 | 19,998 | 10,000 | 7,777,772 | 5,000,000 |
| 21 | 29,998 | 10,000 | 17,777,772 | 10,000,000 |
| 22 | 39,998 | 10,000 | 27,777,772 | 10,000,000 |
| 23 | 49,998 | 10,000 | 77,777,772 | 50,000,000 |
| 24 | 99,998 | 50,000 | 177,777,772 | 100,000,000 |
| 25 | 199,998 | 100,000 | 277,777,772 | 100,000,000 |
| 26 | 299,998 | 100,000 | 777,777,772 | 500,000,000 |
| 27 | 399,998 | 100,000 | 1,777,777,772 | 1,000,000,000 |
| 28 | 499,998 | 100,000 | 2,777,777,772 | 1,000,000,000 |
| 29 | 999,998 | 500,000 | 7,777,777,772 | 5,000,000,000 |
| 30 | 1,999,998 | 1,000,000 | 17,777,777,772 | 10,000,000,000 |
| 31 | 2,999,998 | 1,000,000 | 27,777,777,772 | 10,000,000,000 |

So using negative parts is 9,259 times more efficient than not using them (27,777,777,772 / 2,999,998 = 9,259).

And there are those sevens again; yet another sign from their gods for the mystical Sumerians.

# Appendix S: School Use Today

The Salamis Tablet abacus could be used in education to aid in the teaching of arithmetic, place value, exponents, and scientific notation.

On page 95, Pullan (1968) writes,

> “Much attention has been given to ways in which 'number concepts' develop in a child's mind, and one result of this has been the appearance of a variety of types of 'structural apparatus'. Properly used, they can undoubtedly be of great value, and many teachers have implicit trust in the particular type of apparatus they have chosen. There is, however, a fundamental objection to most of these devices. They give a different picture to each dimension, units, tens, hundreds, etc., whereas it is a basic principle of the Arabic system of notation that the same figures are used in each position. The naughts in such numbers as 10, 100, and 1000 are not intended to give a different appearance but to put the significant figure (in this case 1) in its proper place. They could be omitted if the position of the figure were known ... This difficulty does not arise on the abacus where the same counters 'sometimes stand for more, sometimes for less'; indeed, it is difficult to see what advantages the newly invented types of apparatus have over the old method of counter-casting.”

Other positive features of using The Salamis Tablet in education include:

1. error checking an entry before accumulating, which is more forgiving for the young and/or novice (entering and accumulating are combined on constrained bead abaci like the Soroban, so if an error is made there is no indication and the error is propagated to the answer);

2. positive and negative numbers can be represented and manipulated simultaneously (without using complements as on constrained bead abaci);

3. the methods to accomplish all four arithmetic operations are much more transparent and straightforward than on constrained bead abaci;

4. numbers in near-scientific notation can be represented and manipulated (not possible on constrained bead abaci);

5. all four arithmetic operations, including multiplication as repeated addition and division as repeated subtraction, can be performed without memorizing tables;

6. advanced students can do all the above in base 10 decimal (Roman, Greek, Egyptian), base 12 duodecimal (Roman fractions), or base 60 sexagesimal (Babylonian) numbers (not possible on constrained bead abaci); and

7. much lower cost than either the apparatus discussed by Pullan or constrained bead abaci; i.e., for each abacus just draw some lines on paper and use less than 110 pennies as pebbles to perform all four arithmetic calculations on numbers of the form $a \times b^c$ where b is in the set {10, 12, 60}, 0 <= c < b^(2(1+(b=10))), 1/b <= a < 1, and a has 5(1+(b=10)) significant digits.

## Appendix T: Takeaway

Pebbles (calculos in Latin) are the "bits" used in the Ancients' four function calculator / computer.

The Ancient Computer's normal mode is to work with numbers in what we would call exponential notation. Decimal numbers can have up to 10 significant digits in the coefficient (a fraction < 1 with no leading zeros) and up to 4 significant digits in the exponent (a radix shift). Duodecimal and sexagesimal numbers can have up to 5 significant digits in the coefficient and up to 2 significant digits in the exponent.

Coefficients and exponents can be either positive or negative.

Built-in error checking is included since an addend can be entered and checked before accumulation.

The Ancient Computer is time tested; it or its predecessors have been in use since before 2000 BC.

# Appendix U: Updates

## Frontinus' Duodecimal Abacus

Frontinus' value of "1 1/2 twelfths of a quinaria plus 3/288 plus 2/3 of 1/288 more" indicates a different structure to the Roman duodecimal abacus than described above; i.e.:

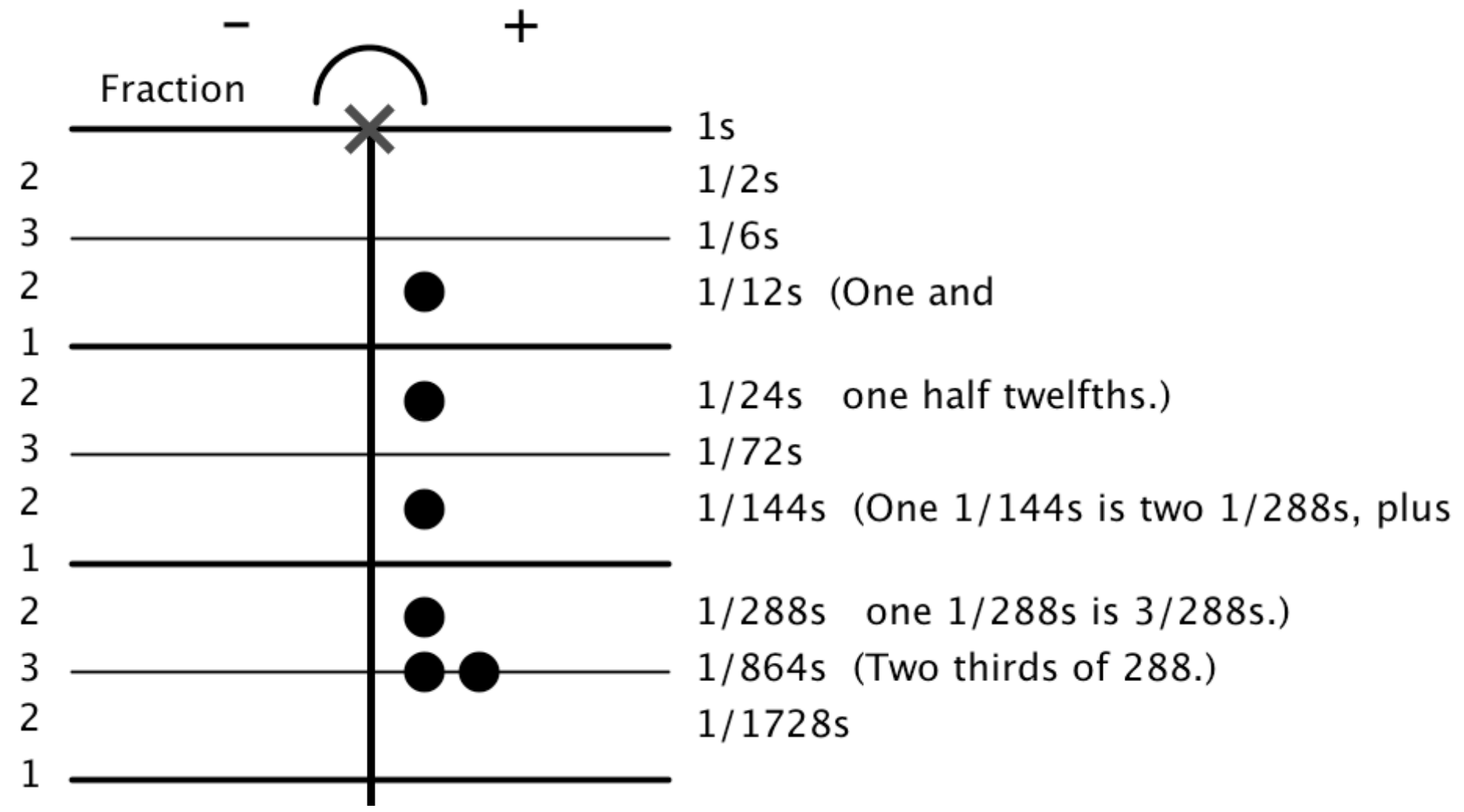

So this is probably the promotion factor structure used by Frontinus. It certainly resolves the lack of knowledge reported on p.388 of Maher and Makowski (2001), "How Frontinus arrived at these particular fractions is not known." Now we know; by reading an answer to a calculation done on a duodecimal mode Salamis Tablet abacus.

This redefinition of the Roman duodecimal line abacus does not change the conclusions in this paper. And it strengthens the argument that both the Romans and Babylonians used The Salamis Table, because to create their duodecimal abacus all the Romans did was change the promotion factors of 5 on the Babylonian abacus to 1.

Changing the Stephenson videos 10.1-10.3 to reflect this new duodecimal abacus is a major effort so will not be done soon. However, the validity of being able to do Frontinus' calculation on an abacus is not diminished, so the existing videos are still good demonstrations.

[1] http://www.cornell.edu
[2] http://www.uml.edu
[3] On pages 17-18, Pullan (1968) writes, "The Latin 'abacus', derived from the earlier Greek word 'abax', meant, simply, 'a flat surface'". The word 'abacus' did not derive from the Hebrew word 'abaq', dust, "and there is little evidence to support a common idea that a table strewn with dust, or sand, was at one time widely used for reckoning." On the other hand, "Sanded tables certainly seem to have been used for the drawing of geometrical figures, ... But it is not so easy to imagine counters being moved easily from place to place on a sandy surface." Menninger, p.301, writes, "The [Romans] called the counting board or table ... the abacus; the Greek work abax means 'round platter' or 'stemless cup' and thus also a 'table without legs.' It is unlikely that this word was derived from the Semitic 'abq' (dust) ..."
[4] Constrained bead abaci have: no addend error checking, no mixed positive and negative numbers, no ability to handle multiple number bases, and no exponents.
[5] On page 89, Pullan (1968) writes, "It is, strictly speaking, a misuse of the word [abacus] to apply it to a bead-frame calculating device."
[6] http://www.materials.ac.uk/awareness/building/index.asp
[7] "In 70AD the Romans laid siege to Jerusalem. The most fearsome weapon in their armory was a massive torsion spring catapult that towered over eight meters high. It could fire boulders weighing twenty-six kilos that destroyed Jerusalem's walls and was so powerful that only the Roman Army's elite tenth legion had the skill to build it. ... Roman eyewitness accounts tell of how devastating the monster catapult was - and how it could backfire, killing their own men as well as the opposition if even the smallest detail was wrong." (Retrieved 8/21/2010 from http://www.materials.ac.uk/awareness/building/romancatapult.asp.)
[8] Watch first minute of http://youtu.be/q_dHpLAPM5I
[9] On page 94, Pullan (1968) writes, " ... it is rare to find ... anyone ... who understands how calculations were made and accounts kept before the introduction of Arabic figures. Too often it is stated, even in authoritative books, that when Roman figures were used there must have been difficulty because of the absence of a sign for zero. This would have been true only if people of the time had been foolish enough to try to do their 'sums' in the same manner that we now use Arabic figures, i.e. by writing on paper. Instead, it was the practice to ... perform the actual calculation with counters on the abacus or counting-board, and to read off the result from the counters as they lay on the board. It could then be recorded on paper in the written notation."
[10] Or use this Line Abacus Worksheet, http://bit.ly/sks23cuLineAbacus.
[11] MacTutor History of Math says that the Akkadians invented the abacus about 2300 BC. (See http://www-history.mcs.st-and.ac.uk/HistTopics/Babylonian_mathematics.html) Ifrah (2000, p. 133) writes, "There is at this time no doubt that the abacus indeed existed in Mesopotamia, and even coexisted with the archaic system of calculi, most probably throughout the third millennium BCE. ... the use of the abacus gave rise to a guild ... the caste of the professional abacists, who no doubt jealously preserved the secrets of their art." Keeping the abacus secret would explain why little documentation of abacus use survives. The maintenance of secrecy about technological advances extends to any State, Ancient or Modern. When the author worked for what was then an advanced computer manufacturer, a license from the State Dept. was needed for every computer shipped outside the U.S.

---

[12] There are: 59 one digit #s; 59 x 59 two digit #s; 59 x 60 x 59 three digit #s; 59 x 60 x 60 x 59 four digit #s; and 59 x 60 x 60 x 60 x 59 five digit #s. Total = 59 ( 1 + 59 ( 1 + 60 + 60^2 + 60^3 ) ) = 59 ( 1 + 59 S ), where 60 S - S = (60 + 60^2 + 60^3 + 60^4 ) - (1 + 60 + 60^2 + 60^3 ), so S = (60^4 - 1) / (60 - 1). Total = 59 ( 1 + 59 (60^4 - 1) / 59 ) = 59 × 60^4. Total without embedded zeros = 59 ( 1 + 59 + 59^2 + 59^3 + 59^4 ) = 59 ( 59^5 - 1 ) / ( 59 - 1 ) = 59 ( 59^5 - 1 ) / 58. Probability that a number has no embedded zeros = ( 59 ( 59^5 - 1 ) / 58 ) / ( 59 x 60^4 ) = ( 59^5 - 1 ) / ( 58 x 60^4 ). Probability that a number has embedded zeros = 1 - ( 59^5 - 1 ) / ( 58 x 60^4 ) ≈ 0.048898070987654 ≈ 5%, or 1 in 20.

[13] Archimedes is often credited with the invention of exponents, due to his treatise, *The Sand Reckoner*. His life spans c. 287 - 212 BC; so he was born after The Salamis Tablet was crafted in marble, c. 300 BC. Was Archimedes' "invention of exponents" inspired by The Salamis Tablet? See http://en.wikipedia.org/wiki/Archimedes, http://www-groups.dcs.st-and.ac.uk/~history/Biographies/Archimedes.html, and http://en.wikipedia.org/wiki/The_Sand_Reckoner.

[14] Smith, G.R. (Smith, 2000), argues that the Hindus originated the sexagesimal system based on their calendar and its inherent astronomical knowledge. Even Ifrah, p.102, writes, "the oldest Iranian civilization" were the Elamites, who spoke a language "some linguists think ... belongs to the Dravidian group (southern India) ... ."

[15] http://en.wikipedia.org/wiki/Egyptian_calendar

[16] Smith, D.E. (1958, p. 160) writes, "That the Egyptians used an abacus is known on the testimony of Herodotus, who says that they 'write their characters and reckon with pebbles, bringing the hand from right to left, while the Greeks go from left to right.' " Also see Appendix E: Egyptian-Babylonian Link.

[17] Retrieved 7/4/2010 from http://en.wikipedia.org/wiki/Square_root_of_2#History

[18] The Old Babylonian period is the time where we get most of our clay tablet artifacts and knowledge of Babylonian mathematics (Melville, 1999). Most of those tablets seem to be results of student work, but frustrating to the historians is the fact that the student's show no work. How did they get their answers?

[19] Image and text retrieved 8/26/2010 from http://philosophy.uchicago.edu/courses/2007-2008-autumn.html and http://www.counton.org/museum/floor4/gallery11/gal11_p2.html respectively.